\documentclass{article}
\usepackage{a4wide}
\usepackage{natbib}

\usepackage[latin1]{inputenc}

\usepackage{amsmath}
\usepackage{amssymb}
\usepackage{graphicx}
\usepackage{url}
\usepackage{xcolor}
\usepackage[draft,silent]{fixme}

\usepackage{lscape}
\usepackage[T1]{fontenc}
\usepackage[sc]{mathpazo}
\def \1{\mathbf{1}}
\def\inproba{\stackrel{p}{\longrightarrow}}
\def\inlaw{\stackrel{{d}}{\longrightarrow}}
\newcommand{\inv}[1]{\frac{1}{#1}}
\renewcommand{\O}{\mathcal O}
\newcommand{\floor}[1]{\lfloor #1 \rfloor}

\def \P{\mathbb{P}}
\def \PE{\mathbb{E}}
\def\rmd{\textrm{d}}

\DeclareMathOperator*{\cov}{Cov}

\DeclareMathOperator*{\var}{Var}
\DeclareMathOperator*{\argmax}{\operatorname{argmax}}

\newtheorem{theorem}{Theorem}

\newtheorem{coro}{Corollary}

\graphicspath{{fig/}}

\title{Homogeneity and change-point detection tests for
  multivariate data using rank statistics}
\author{
Alexandre Lung-Yut-Fong, Céline Lévy-Leduc, Olivier Cappé
}
\begin{document}
\maketitle              % typeset the title of the contribution

\begin{abstract}
  Detecting and locating changes in highly multivariate data is a major concern in several current
  statistical applications.
  In this context, the first contribution of the paper is a novel
  non-parametric two-sample homogeneity test for multivariate data based on the well-known Wilcoxon
  rank statistic. The proposed two-sample homogeneity test statistic can be extended to deal with
  ordinal or censored data as well as to test for the homogeneity of more than two samples. We also
  provide a detailed analysis of the power of the proposed test statistic (in the two sample case)
  against asymptotic local shift alternatives.
  The second contribution of the paper concerns the use of the proposed test statistic to perform
  retrospective change-point detection. It is first shown that the approach is computationally
  feasible even when looking for a large number of change-points thanks to the use of dynamic
  programming. Computable asymptotic $p$-values for the test are available in the case where a
  single potential change-point is to be detected.
  The proposed approach is particularly recommendable in situations where the correlations between
  the coordinates of the data are moderate, the marginal distributions are not well modelled by
  usual parametric assumptions (e.g., in the presence of outliers) and when faced with highly variable
  change patterns, for instance, if the potential changes only affect subsets of the coordinates of
  the data.
\end{abstract}

%\keywords{change-point detection; homogeneity test; Kruskal-Wallis test; Mann-Whitney/Wilcoxon test; multivariate data; rank statistics}

\section{Introduction}
\label{sec:introduction}

Detection and location of distributional changes in data is a major
statistical challenge that arises in many different contexts. This very general
concern can be particularised to more specific tasks such as segmentation,
novelty detection or significance tests. In this contribution, we focus on two
types of problems: \emph{homogeneity testing}, where the statistician is
presented with pre-specified groupings of the data that are believed to be
comparable, and \emph{change-point detection}, in which a series --most often,
a time series-- is to be segmented into homogeneous contiguous regions. These
two tasks are obviously related but the latter is more challenging as the
appropriate groupings of the data are unknown, although one does have the
strong prior assumption that homogeneous regions of the data are
contiguous. Homogeneity testing and/or change-point detection are instrumental
in applications that range from the surveillance of industrial processes
\citep{Basseville:1993}, to computer
security \citep{Tartakovsky:2006,Levy:Roueff:2009}, processing of audiovisual
data \citep{desobry:davy:doncarli:2005}, financial and econometric modelling
\citep{bai:perron:2003,talih:hengartner:2005}, health monitoring
\citep{brodsky:2000}, or  bioinformatics
\citep{picard:robin:lavielle:vaisse:daudin:2005,vert:bleakley:2010}.

In light of the important available literature on change-point detection it is
important to make two additional distinctions. First, in many cases the data to
be analysed can be assumed to present some form of global reproducibility and to
include several instances of actual changes. In this case, it seems reasonable
to fit a model to the data to profit from the available statistical
information regarding various relevant aspects of the problem such as the
distribution of the data in the absence of change, the typical change-point
patterns, etc. In such situations, very convincing results have been
demonstrated using Bayesian approaches due to the existence of efficient
computational methods to explore the posterior distribution, even when using very
flexible models \citep{barry:hartigan:1992,fearnhead:2006}. In contrast,
in this contribution, we consider scenarios in which the data are either scarce
or very variable or where potential changes occur somewhat infrequently. In
this alternative context, the goal is to develop approaches that make as few
assumptions as possible regarding the underlying distribution of the data or the
nature of the changes and that do not rely on the observation of actual
change patterns. The second important distinction is that many works in the time
series literature consider the \emph{online} change-point detection
framework \citep{siegmund:1947}
in which the data have to be processed on-the-fly or with minimal delay
using for instance the CUSUM algorithm initially proposed by
\cite{page:1954}.
%\fixme{Why is this called Page-Hinkley?}
In the following, we consider the opposite
situation, sometimes referred to as \emph{retrospective} analysis, in which all
the data to be tested have been recorded and are available for analysis.

In this context, the first contribution of this work consists in novel homogeneity tests for
dealing with possibly high-dimensional multivariate observations. We focus on situations where the potential changes are believed to have a strong impact on the mean but where the distribution of the data is otherwise mostly unknown. For scalar observations, there
are well-known robust solutions for testing homogeneity in this context such as the Wilcoxon/Mann--Whitney or
Kruskal-Wallis procedures \citep{lehmann:1975} to be further discussed below. For multivariate
observations, the situation is far more challenging as one would like to achieve robustness with
respect both to the form of the marginal distributions and to the existence of correlations (or
other form of dependence) between coordinates. The latter aspect has been addressed in a series of
works by
\cite{mottonen:oja:tienari:1997,hettmansperger:mottonen:oja:1998,oja:1999,topchii:tyurin:oja:2003}
who studied multivariate extensions of sign and rank tests. These tests are affine invariant in the
sense that they behave similarly to Hotelling's $T^2$ test (see Section~\ref{sec:power} below)
---which is optimal for Gaussian distributions-- for general classes of multivariate distributions
having ellipsoidal contours (e.g., for multivariate $t$ distributions). Another promising approach
investigated by several recent works consists in using kernel-based methods
\citep{desobry:davy:doncarli:2005,gretton:2006:nips,harchaoui:bach:moulines:2008}. In our
experience however, these methods that can achieve impressive results for moderately
multidimensional data or in specific situations (\textit{e.g}., if the data lie on a
low-dimensional manifold) lack robustness when moving to larger dimensions. In particular, as
illustrated in Section~\ref{sec:simulToy} below, kernel-based methods are not robust with respect
to the presence of contaminating noise and to the fact that the changes to be detected may only
affect a subset of the components of the high-dimensional data. The latter scenario is of very
important practical significance in applications where the data to be analysed consist in
exhaustive recordings of complex situations that are only partly affected by possible changes (see
\cite{Levy:Roueff:2009} for an example regarding the detection of computer attacks). The method
proposed in this work is based on a combination of marginal rank statistics, following the
pioneering idea of \cite{wei:lachin:1984}. Compared to the latter, our contribution is twofold:
first, we show how to correct the bias that appears in the test statistic proposed by
\cite{wei:lachin:1984} whenever the two samples are not balanced in size; we then show how to
extend this idea to the case of more than two groups. The numerical simulation presented in
Section~\ref{sec:Eval1} confirms that the proposed test statistic is significantly more robust than
kernel-based methods or approaches based on least-squares or Hotelling's $T^2$ statistics. These
empirical observations are also supported by the results detailed in Section~\ref{sec:power}
regarding the power of the proposed test statistic (in the two sample case) with respect to asymptotic
local shift alternatives. In particular, although the proposed test is not strictly affine
invariant\footnote{It is highly unlikely that there exist statistics which are both invariant under
  monotonic transformations of the marginals --as the proposed method is-- and affine invariant.} it
compares favourably to Hotelling's $T^2$ statistic for important classes of distributions under the
assumption that the condition number of the correlation matrix of the data is not too large.

We then consider the use of the proposed approach for change-point detection by
optimising the test statistic over the --now considered unknown-- positions of
the segment boundaries. Although simple this idea raises two type of
difficulties. The first one is computational as the resulting optimisation
task is combinatorial and cannot be solved by brute force enumeration when
there is more than one change-point (that is, two segments). In the
literature this issue has been previously tackled either using dynamic
programming \citep{bai:perron:2003,Harchaoui:Cappe:2007} or more
recently using
Lasso-type penalties \citep{harchaoui:levy:2010,vert:bleakley:2010}. We show
that the generic dynamic programming strategy % as described by \cite{Bellman:1961,kay:1993}
is applicable to the proposed test statistic
making it practically suitable for retrospective detection of multiple
change-points. The second difficulty is statistical as the optimisation
with respect to the change-point locations modifies the distribution of the
test statistics. Thus, the design of quantitative criterions for assessing the
significance of the test is a challenging problem in this context. This
issue has been considered before, mostly in the case of a single
change-point, for various test statistics
%\fixme{Is thereother general refs. on changepoints?}
\citep{Csorgo:1997,chen:gupta:2000}. In many cases the asymptotic distribution of the test remains
hard to characterise and must be calibrated using Monte Carlo simulations.
We show that for a simple modification of the proposed
rank-based statistic, one can indeed obtain computable asymptotic $p$-values
that can be used to assess the significance of the test when looking for a
single change-point.

The paper is organised as follows. Section~\ref{sec:TestingHomo} is devoted to
homogeneity testing, starting first with the two-sample case and then
considering the more general situation where several predefined groups are
available. The end of Section~\ref{sec:TestingHomo} is dedicated
 to the study of the asymptotic behaviour of the two-sample
 homogeneity test under local shift alternatives.
In Section~\ref{sec:ChangePoint}, the proposed test statistic is modified to provide a method for detecting and locating change-points, with the computation of $p$-values (in the single change-point case) being discussed in Section~\ref{sec:CP-onechange}.
The results of numerical experiments carried out both on simulated
and on real data are then reported in Section~\ref{sec:Eval1}.
% The appendices contain the proofs of the theorems stated in the text.

\section{Testing for Homogeneity}
\label{sec:TestingHomo}
We first tackle in this Section the so-called two-sample problem, that is testing the homogeneity
between two partitions of data.
The proposed test statistic is then extended in Section \ref{sec:test_homogeneity_multigroups} to deal with more than two groups of data.
\subsection{Two-sample homogeneity test}
\label{sec:test_twosample}
% In this section, we first consider testing homogeneity based on rank statistics in a multivariate
% setting.
Consider $n$ $K$-dimensional multivariate observations $(\mathbf{X}_{1}, \dots,
\mathbf{X}_{n})$ and denote by $X_{i,k}$ the $k$th coordinate of $X_i$, such
that $\mathbf{X}_{i}=(X_{i,1},\dots,X_{i,K})'$, where the prime is used to
denote transposition. 
%We first consider testing homogeneity of two subsections
%of the data before extending the proposed test statistic to deal with more than
%one potential change in Section~\ref{sec:test_homogeneity_multigroups}. 
%To this aim,
We consider
the classical statistic test framework with the null (or baseline) 
hypothesis, $(H_{0})$: ``$(\mathbf{X}_{1},
\dots, \mathbf{X}_{n})$ are identically distributed random vectors'', and
the alternative hypothesis, $(H_{1})$: ``$(\mathbf{X}_{1}, \dots, \mathbf{X}_{n_1})$ are distributed under $\Bbb P_{1}$ and
$(\mathbf{X}_{n_1+1}, \dots, \mathbf{X}_{n})$ under $\Bbb P_{2}$, with $\Bbb P_{1} \neq \Bbb P_{2}$''. In this setting, the potential change point $n_1$ is assumed to be given but the data distributions are fully unspecified both under $(H_{0})$ and $(H_{1})$. The proposed test statistic extends the well-known Wilcoxon/Mann--Whitney rank-based criterion to multivariate data by considering the asymptotic joint behaviour of the rank statistics that can be computed from each coordinate of the observations. For $k$ in $\{1,\dots,K\}$, define the vector-valued statistic $\mathbf{U}_{n}(n_1)=(U_{n,1}(n_1),\dotsc,U_{n,K}(n_1))'$ by
\begin{equation}\label{eq:Ukn}
  U_{n,k}(n_1) =
 \inv{\sqrt{n n_1 (n-n_1)}} \sum_{i=1}^{n_{1}} \sum_{j=n_{1}+1}^{n} 
   \left\{\1(X_{i,k}\le X_{j,k})-\1(X_{j,k}\le X_{i,k})\right\} \;.
\end{equation}
Although the form of the statistic given above is more appropriate for mathematical analysis as well as for discussing possible generalisations of the approach (see Section~\ref{sec:ext} below), it is important to realise that $U_{n,k}(n_1)$ is related to the classical Wilcoxon/Mann--Whitney statistic computed from the series $X_{1,k}, \dots, X_{n,k}$. Assuming that there are no ties in the data, let $R_j^{(k)}$ denote the rank of $X_{j,k}$ among
$(X_{1,k},\dots,X_{n,k})$, that is, 
$R_j^{(k)}=\sum_{i=1}^n \1(X_{i,k}\leq X_{j,k})$. Noticing that $\sum_{j=1}^n R_j^{(k)} = n (n+1)/2$, it is then easily verified that $U_{n,k}(n_1)$ can be
equivalently defined as
\begin{equation}\label{eq:ranks}
 U_{n,k}(n_1)=\frac{2} {\sqrt{n n_1(n-n_1)}}\sum_{i=1}^{n_{1}}\left(\frac{n+1}{2}-R_i^{(k)}\right)
=\frac{2} {\sqrt{n n_1(n-n_1)}} \sum_{j=n_{1}+1}^{n} \left(R_j^{(k)}-\frac{n+1}{2}\right)\;.
\end{equation}
This alternative form of $U_{n,k}(n_1)$ is more appropriate for
computational purposes as discussed in Section~\ref{sec:implDetails}
below. For convenience, we denote by $\hat{F}_{n,k}(t) = n^{-1}\sum_{j=1}^n
\1(X_{j,k}\leq t)$ the empirical cumulative distribution function
(c.d.f.\ in short) of the $k$th coordinate, such that
$\hat{F}_{n,k}(X_{i,k}) = R_i^{(k)}/n$. Let $\hat\Sigma_{n}$ denote the $K$-dimensional empirical covariance matrix defined by
\begin{equation}\label{eq:Sigma_emp}
 \hat\Sigma_{n, kk'} =
  \frac4n\sum_{i=1}^n\{\hat{F}_{n,k}(X_{i,k})-1/2\}\{\hat{F}_{n,k'}(X_{i,k'})-1/2\},\;
  1\leq k,k'\leq K\; .
\end{equation}
The test statistic that we propose for assessing the presence of a potential change in $n_1$ is defined as
\begin{equation}\label{eq:StatWL}
S_{n}(n_1) = \mathbf{U}_{n}(n_1)' \hat{\Sigma}_n^{-1} \mathbf{U}_{n}(n_1)\;.
\end{equation}
Theorem~\ref{th:WL-theorem} below (proved in Appendix~\ref{sec:Appendix}) gives the
limiting behaviour of the test statistic $S_{n}(n_1)$ under the null hypothesis.

\begin{theorem}
  \label{th:WL-theorem}
Let $\mathbf{X}_{1},\dotsc,\mathbf{X}_{n_{1}},\mathbf{X}_{n_{1}+1},\dotsc,\mathbf{X}_{n}$
be $\Bbb R^{K}$-valued i.i.d.\
random vectors, such that for all $k$ in $\{1,\dotsc,K\}$, 
the cumulative distribution function $F_{k}$ of $X_{1,k}$ 
 is a continuous function. Assume that $n_{1}/n \to
 t_1$ in (0,1) as $n$ tends to infinity, and that the $K\times K$ covariance matrix $\Sigma$
 defined by
\begin{equation}
  \label{eq:CovMatrix}
  \Sigma_{k k'} = 4 \cov\left(F_{k}(X_{1,k});F_{k'}(X_{1,k'})\right),\;
  1\leq k,k'\leq K \; ,
\end{equation}
is positive definite. Then, the test statistic $S_{n}(n_1)$ defined in~\eqref{eq:StatWL} converges in
distribution to a $\chi^{2}$ distribution with $K$ degrees of freedom.
\end{theorem}

Theorem~\ref{th:WL-theorem} shows that the proposed test is well normalised
with respect to the dimension $K$, the length $n$ of the data and the
postulated change-point location $n_1$. It is asymptotically distribution-free
in the sense that its limiting behaviour under $(H_0)$ does not depend on the
distribution of the data. By construction, it is also invariant under any
monotonic transformation of the coordinates of $\mathbf{X}_i$.

The matrix $\Sigma$, which corresponds to the asymptotic covariance matrix of
the vector $\mathbf{U}_{n}(n_1)$ is equal, up to a multiplicative constant, to
the Spearman correlation matrix of $X_i$
\citep{lehmann:1975,vandervaart:1998}. This is a well-known robust measure of
dependence that appears in particular when using copula models. A sufficient
condition for ensuring the invertibility of $\Sigma$ is thus that no linear
combination of the $F_k(X_{1,k})$'s should be almost surely equal to a
constant, which is arguably a very weak condition. It is easily checked that
the diagonal elements of $\Sigma$ are all equal to $1/3$ 
% \footnote{It appears
  % that the diagonal elements of $\hat\Sigma_{n}$ converge very rapidly to this
  % value and we did not observe any significant improvement when trying to take
  % into account this fact when estimating $\Sigma$.}
and that $\Sigma_{k\ell} =
\Sigma_{\ell k} = 0$ whenever the $k$-th and $\ell$-th coordinates of
$\mathbf{X}_i$ are independent. It appears, in practice, that the diagonal elements of $\hat\Sigma_{n}$ converge very rapidly to 1/3
  value and we did not observe any significant improvement when trying to take
  into account this fact when estimating $\Sigma$.

Theorem~\ref{th:WL-theorem} defines the asymptotic false alarm rate associated
with the test statistic $S_n(n_1)$. The test is consistent (\textit{i.e.},
its power tends to 1) for all alternatives such that the condition ensuring the
consistency of the standard Wilcoxon/Mann--Whitney two-sample test holds true
for at least one coordinate. More formally, the proposed test is consistent when
there exists $k$ in $\{1,\dots,K\}$ such that $\P(X_{k,1}\leq X_{k,n})\neq
1/2$. In the scalar case, this condition is known to hold for general classes of changes such as
shift (change-in-the-mean) models or scale (multiplicative) change for positive
variables. We defer to Section~\ref{th:power} a more detailed investigation of the asymptotic power of the test in the case of multivariate shift alternatives.

\begin{figure}[htb] \centering
  \begin{tabular}{ccc}
    \includegraphics[width=.31\textwidth]{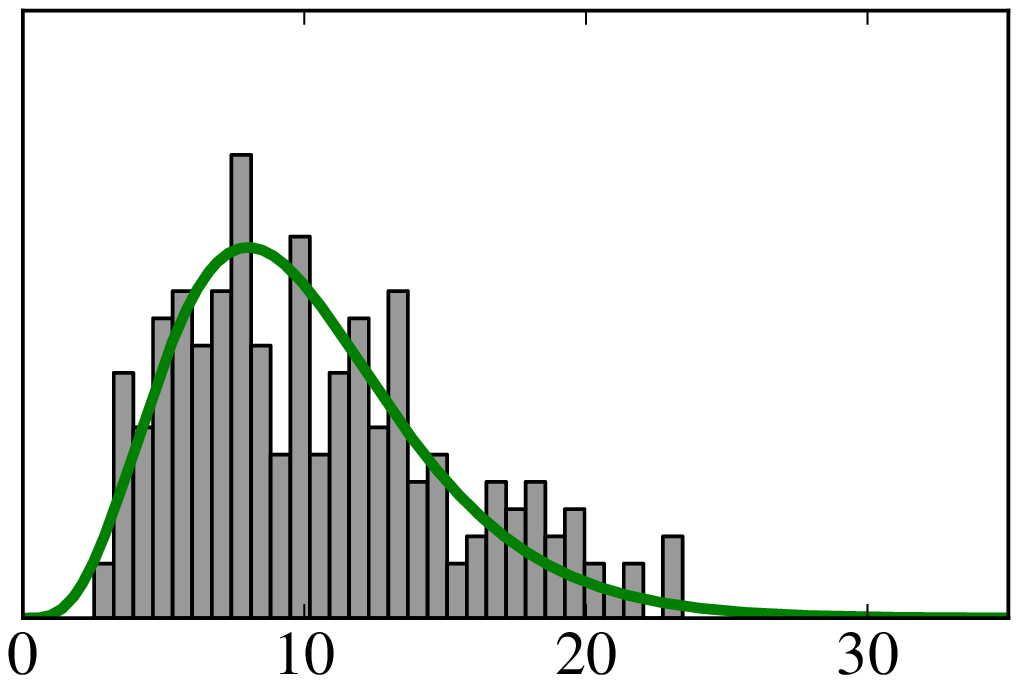}&
    \includegraphics[width=.31\textwidth]{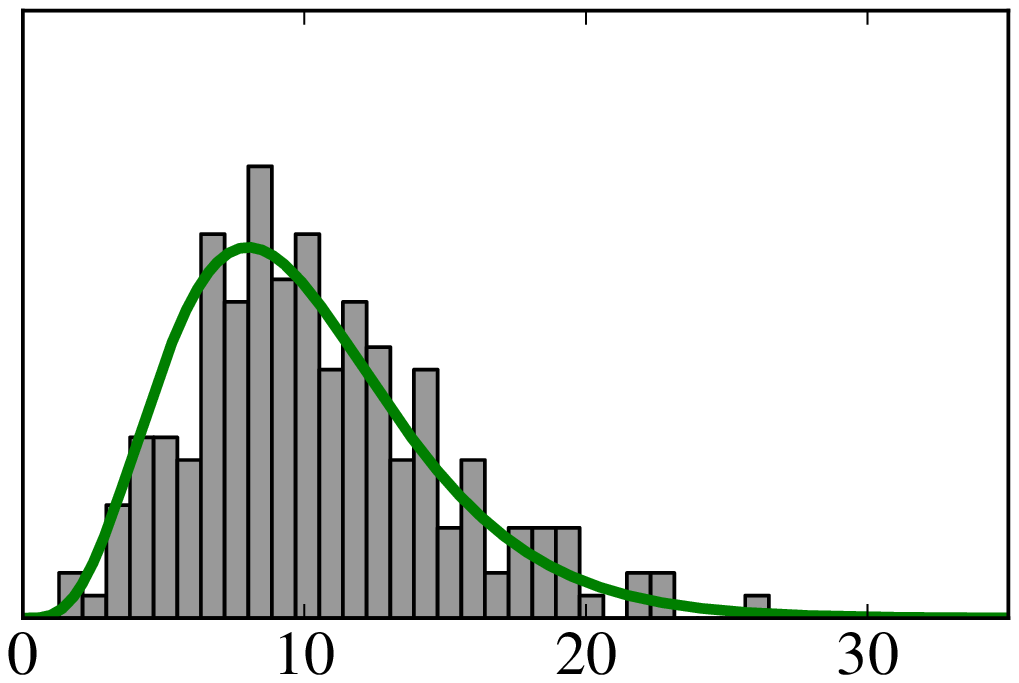}&
    \includegraphics[width=.31\textwidth]{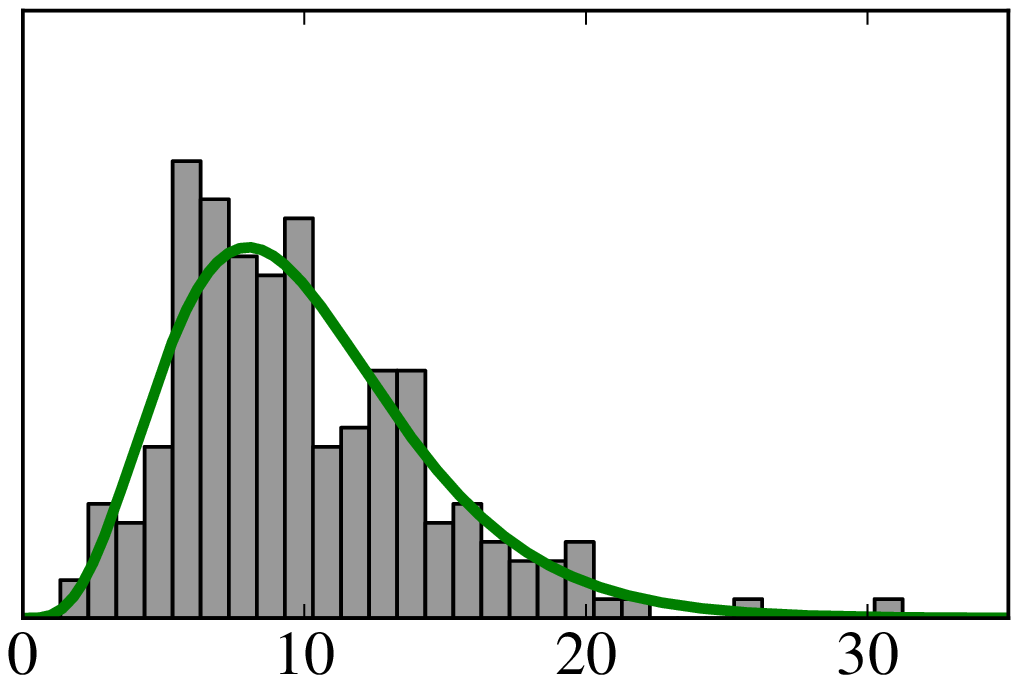}\\
    \includegraphics[width=.31\textwidth]{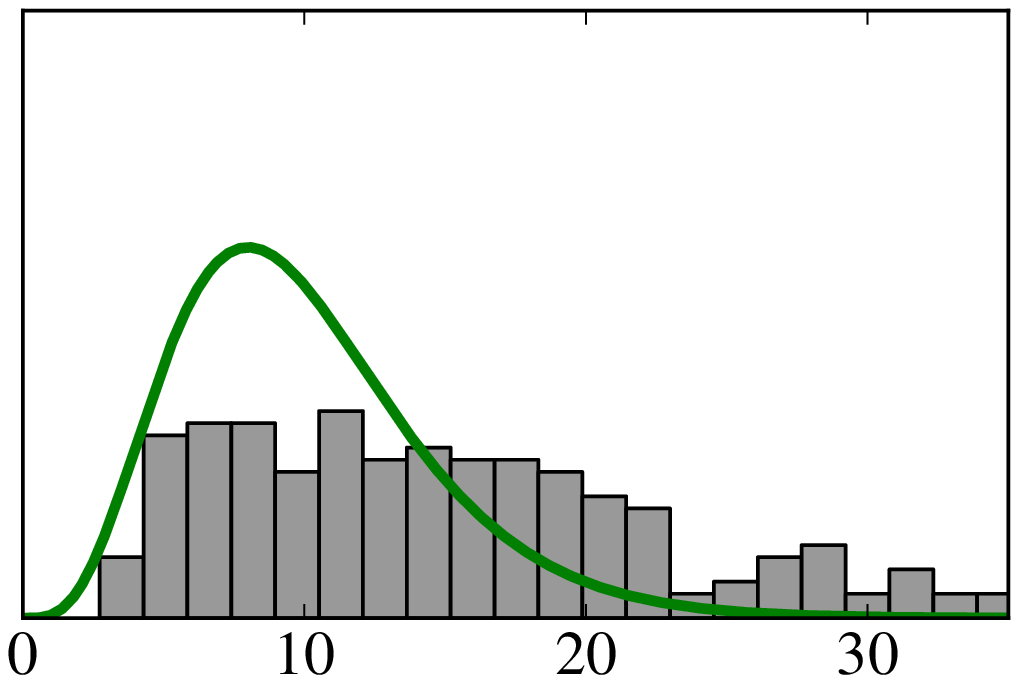}&
    \includegraphics[width=.31\textwidth]{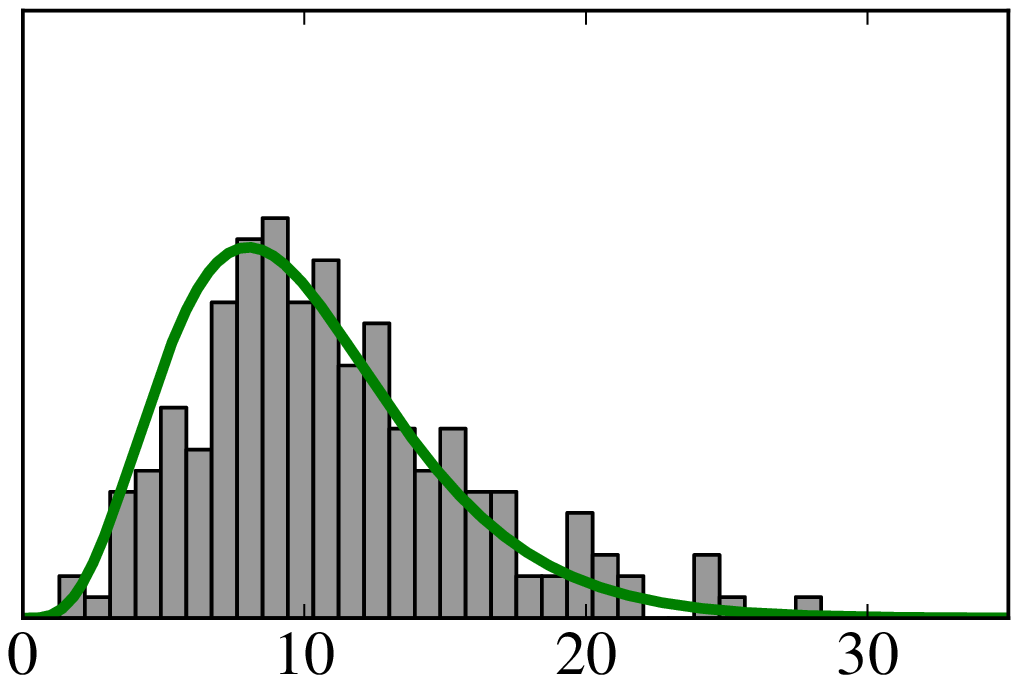}&
    \includegraphics[width=.31\textwidth]{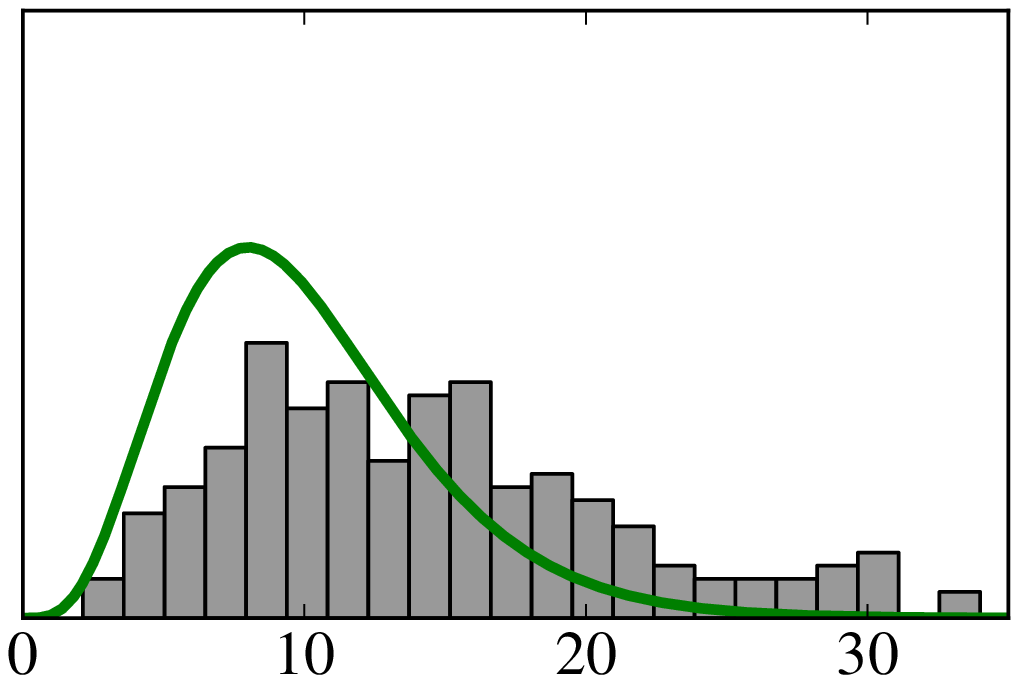}\\
    $n_1/n=1/8$ & $n_1/n=1/2$ & $n_1/n=7/8$
  \end{tabular}
  \caption{Histograms of the statistics $S_n(n_1)$ (top) and Wei and
    Lachin's (bottom) compared to the
    $\chi^2_K$ p.d.f., as a function of the ratio $n_1/n$. Data
    corresponds to $n=200$ samples of a ten-dimensional standard Gaussian distribution.}
  \label{fig:Comparison_WL_WLmaison}
\end{figure}

As Theorem~\ref{th:WL-theorem} is an asymptotic result, we have carried out
Monte Carlo simulations to asses the accuracy of the approximation
for finite sample sizes. Using data with
independent coordinates\footnote{Note that by construction, the test statistic is then fully
invariant with respect to the precise distribution used for the Monte Carlo
simulations.} 
we found that the distribution of $S_n(n_1)$ defined in
(\ref{eq:StatWL}) can be considered close enough
to the limiting distribution, as measured by the Kolmogorov-Smirnov
test at level 1\%, when $n$ is at least 8 times larger than $K$.
 For instance, for $K=20$, a value of $n=210$ was sufficient; $K=100$ required
$n=840$ samples, etc. The empirical density of the test statistics is illustrated in the upper part of Figure~\ref{fig:Comparison_WL_WLmaison} when $K=10$ and $n=200$.

The pioneering work of \cite{wei:lachin:1984} describes a result analogous to
that of Theorem~\ref{th:WL-theorem} in the case of possibly upper-censored
data. However, the proof technique used by \cite{wei:lachin:1984} relies on a
different interpretation of $\Sigma$ which is used to derive a weighting
matrix that is not equal to $\hat\Sigma_{n}$ as defined
in~\eqref{eq:Sigma_emp}. In contrast, our proof (see Appendix~\ref{sec:Appendix})
is based on a standard argument for U-statistics (the
Hoeffding decomposition) that directly returns an expression of $\Sigma$ in
terms of covariances, for which usual estimators, such as $\hat\Sigma_{n}$, may
be used. The test statistics of \cite{wei:lachin:1984} thus differs from
$S_n(n_1)$ and turns out to be biased in cases where $n_1 \neq n/2$, \textit{i.e.}, when
the change does not occur in the centre of the observation frame, as shown by
the bottom part of Figure~\ref{fig:Comparison_WL_WLmaison}. This bias becomes
problematic when the potential change location $n_1$ is unknown because the values
of $S_n(n_1)$ for different values of $n_1$ cannot be validly compared.

\subsubsection{Implementation Issues}
\label{sec:implDetails}

As noted above, the vector $(U_{n,k}(n_1))_{1\leq k\leq K}$ should be computed from the marginal rank statistics in the form given in (\ref{eq:ranks}). $\hat\Sigma_{n}$ also is a simple function of those marginal ranks. Thus, $(U_{n,k}(n_1))_{1\leq k\leq K}$ can be computed in  $\O(K n\log(n))$ operations using a sort for computing the ranks, as the average numerical complexity of usual sorting algorithms is of the order of $n\log(n)$ operations. The computation of
$\hat{\Sigma}_n$ then requires $\O(K^{2}n)$ operations and its inversion
$\O(K^{3})$ operations. Note that if the test statistic needs to be
recomputed at a neighbouring index, say $n_1+1$, neither the ranks nor
$\hat{\Sigma}_n$ and its inverse need to be recomputed. Hence the number of
additional operations required to compute $\hat{S}_n(n_1+1)$ is indeed very limited.

%\fxnote{Dager sign also used for footnotes}
In some situations, it may happen that the empirical estimate $\hat{\Sigma}_n$
becomes ill-conditioned rendering its inversion numerically
unstable. \cite{wei:lachin:1984} suggested to circumvent the problem by adding
some small positive value to the diagonal elements of $\hat{\Sigma}_n$. It is
important to realise however that a particular case where $\Sigma$ itself can
be ill-conditioned is when coordinates of $\mathbf{X}_{1}$ are strongly
dependent. In the limiting case where coordinates of $\mathbf{X}_{1}$, say
two of them for illustration, are duplicated, $\Sigma$ becomes a matrix of rank
$K-1$. In such a case, the correct statistic is
obtained by simply discarding one of the coordinates that are duplicated. Hence,
to regularise $\hat{\Sigma}_n$ in such cases, we suggest inverting it using its
Moore-Penrose pseudo inverse: if $\hat{\Sigma}_n = USU'$ denotes the singular
value decomposition of $\hat{\Sigma}_n$, with $S =
\operatorname{diag}(s_{1},\dots,s_{K})$ being the diagonal matrix of
eigenvalues of $\Sigma$, then the pseudo inverse $\Sigma_n^{\dagger}$ is
defined as $U'\operatorname{diag}(s_{1}^{\dagger},\dots,s_{K}^{\dagger}) U$
where $s_{i}^{\dagger} = s_{i}^{-1} \1(s_{i} >\epsilon)$ and $\epsilon$ is a
fixed positive threshold. Instead of relying on the asymptotic result of
Theorem~\ref{th:WL-theorem}, it is suggested to compare $S_{n}(n_1)$ to the
quantiles of the $\chi^{2}_{K'}$ distribution, where $K'$ is the number of
non-null values among the $s_{i}^\dagger$'s. As already mentioned however, some terms
of $\hat{\Sigma}_n$ appear to converge very rapidly and the matrix is only very
rarely ill-conditioned, even when $n$ is only slightly larger than $K$. On the
other hand, the regularised variant described above was found to be
effective for dealing with
signals whose coordinates can be extremely dependent, \textit{e.g.}, if there is a quasi-deterministic relationship between two coordinates.

% Especially when dealing with very large dimensional signals, \textit{i.e.}, when $K$ is of the
% order of $n$, a different option might be to take profit of the
% observation that the diagonal elements of $\Sigma$ are equal to $1/3$ (assuming
% continuous c.d.f.\'s) to use a shrinkage estimate of $\hat{\Sigma}_n$, with $1/3
% \mathbf{Id}_K$ as the default shrinkage target \cite{ledoit:wolf:2004}. This approach,
% in which the shrinkage level is optimised w.r.t. to some covariance matching criterion is
% preferable to the ``blind'' diagonal loading regularization proposed in \cite{wei:lachin:1984}.
% For the type of dimensions considered in the examples of Section~\ref{sec:Eval1} ($K$
% equal to a few tens and $n$ to a few hundreds), this modification was not found
% to be beneficial.

\subsubsection{Discrete, missing or censored data}
\label{sec:ext}
Theorem~\ref{th:WL-theorem} requires the continuity of the c.d.f.\ $F_{k}$ of each
coordinate; hence it is not directly applicable, for instance, to discrete variables. In such cases
however, Theorem~\ref{th:WL-theorem} is still valid upon redefining $\Sigma$ as
\begin{equation}
  \label{eq:CovMatrixBis}
  \Sigma_{kk'}  = \PE\left[\{F_{k}(X_{1,k}^{-})+F_{k}(X_{1,k})-1\}
\{F_{k'}(X_{1,k'}^{-})+F_{k'}(X_{1,k'})-1\}\right]\; ,
\end{equation}
where $F_{k}(x^-)$ denotes the left-limit of the c.d.f.\ in $x$. In
this case, (\ref{eq:ranks}) has to be replaced by 
$$
U_{n,k}(n_1) 
=\frac{2}{\sqrt{n n_1 (n-n_1)}}\sum_{j=n_{1}+1}^{n}
\left\{R_j^{(k)}-\frac{n+\sum_{i=1}^{n}\1(X_{i,k}=X_{j,k})}{2}\right\}\;.
$$
% for details, see Appendix~\ref{proof:eq:ranks}.
%\fxnote{Removed pointer to Appendix.}

Another useful extension concerns the case of censored or missing data that can
be dealt with in great generality by introducing lower $\underline{X}_{i,k}$
and upper $\overline{X}_{i,k}$ censoring values such that $\underline{X}_{i,k}
\leq X_{i,k} \leq \overline{X}_{i,k}$, where a strict inequality indicates
censoring (for missing values, simply set $\underline{X}_{i,k}=-\infty$ and
$\overline{X}_{i,k}=+\infty$). In this case, we define a modified statistics from the censoring bounds $\underline{X}_{i,k}$
and $\overline{X}_{i,k}$ by
\begin{equation}
  \label{eq:Ukn-cens}
  U_{n,k}(n_1) =
 \inv{\sqrt{n n_1 (n-n_1)}} \sum_{i=1}^{n_{1}} \sum_{j=n_{1}+1}^{n} 
   \left\{\1(\overline{X}_{i,k}\le \underline{X}_{j,k})-\1(\overline{X}_{j,k}\le \underline{X}_{i,k})\right\} \;.
\end{equation}
It can be shown, adapting the arguments of~\cite{lung:levyledu:cappe:2011} who studied the use of the scalar statistic~\eqref{eq:Ukn-cens} for change-point detection, that Theorem~\ref{th:WL-theorem} then holds with
\begin{equation}
  \label{eq:CovMatrix-cens}
  \Sigma_{kk'}  = \PE\left[\{\overline{F}_{k}(\underline{X}_{1,k})+\underline{F}_{k}(\overline{X}_{1,k}^{\,-})-1\}
\{\overline{F}_{k'}(\underline{X}_{1,k'})+\underline{F}_{k'}(\overline{X}_{1,k'}^{\,-})-1\}\right]\; ,
\end{equation}
where $\overline{F}_{k}$ and $\underline{F}_{k}$ denote the c.d.f.'s
of $\overline{X}_{1,k}$ and $\underline{X}_{1,k}$, respectively.
%\textcolor{red}{and by changing in the definition of $S_n(n_1)$ (see
%  Equation (\ref{eq:StatWL})) $\hat{\Sigma}_n$ defined in
%  (\ref{eq:Sigma_emp}) by the estimator of the new matrix $\Sigma$
%  defined in (\ref{eq:CovMatrix-cens})}.

\subsection{Testing homogeneity within several groups of data}
\label{sec:test_homogeneity_multigroups}

% The change-point estimators that we propose are obtained as
% maximizers of a testing procedure used for testing the equality
% in distribution of several groups as described below.

In this section, the procedure presented so far is extended to deal with more
than two groups of multivariate data. The resulting test statistic is again
based on a proper combination of marginal statistics involved in the
Kruskal-Wallis procedure that generalises the classical Wilcoxon-rank test
when there are more than two groups of data.

Consider the null hypothesis that
$L$ given groups,
$\mathbf{X}_{1},\dots,\mathbf{X}_{n_1}$;
$\mathbf{X}_{n_1+1},\dots,\mathbf{X}_{n_2}$;
\dots; $\mathbf{X}_{n_{L-1}+1},\dots,\mathbf{X}_{n_{L}}$, share
the same distribution, where we shall use the convention that $n_0=0$ and
$n_{L}=n$. 
For $j$ in $\{1,\dots,n\}$ and $k$ in $\{1,\dots,K\}$, denote as
previously by $R_j^{(k)}$ the rank of $X_{j,k}$ among
$(X_{1,k},\dots,X_{n,k})$ that is,
$R_j^{(k)}=\sum_{i=1}^n \1_{\{X_{i,\ell}\leq X_{j,\ell}\}}$. For $\ell$ in $\{0,\dots,L-1\}$, define the average rank in group $\ell$ for the $k$th coordinate
by 
$
\bar{R}_{\ell}^{(k)}=(n_{\ell+1}-n_{\ell})^{-1}\sum_{j=n_{\ell}+1}^{n_{\ell+1}}
R_j^{(k)}.
$
Consider the following test statistic:
\begin{equation}\label{eq:KW_mult}
T(n_1,\dots,n_{L-1})=\frac{4}{n^2}\sum_{\ell=0}^{L-1}(n_{\ell+1}-n_\ell)\mathbf{\bar{R}}_\ell'\;\hat{\Sigma}_n^{-1}\;
\mathbf{\bar{R}}_{\ell}\;,
\end{equation}
where the vector $\mathbf{\bar{R}}_\ell$ is defined as
$
\mathbf{\bar{R}}_\ell=(\bar{R}_{\ell}^{(1)}-(n+1)/2,\dots,\bar{R}_{\ell}^{(K)}-(n+1)/2)',
$ and $\hat{\Sigma}_n$ is again the matrix defined in~(\ref{eq:Sigma_emp}).
% the $K\times K$ matrix 
% whose $(k,k')$--th element is defined by
% \begin{equation}\label{eq:Sigma_emp_rank}
%  \hat\Sigma_{n, k k'} =
%   \frac{1}{n^2}\sum_{i=1}^n\{R_i^{(k)}-n/2\}\{R_i^{(k')}-n/2\}\;.
% \end{equation}
% Note that $ \hat\Sigma_{n, k k'}$ can be rewritten as
% \begin{equation}\label{eq:Sigma_emp_mkw}
%  \hat\Sigma_{n, k k'} =
%   \frac1n\sum_{i=1}^n\{\hat{F}_{n,k}(X_{i,k})-1/2\}\{\hat{F}_{n,k'}(X_{i,k'})-1/2\}\;,
% \end{equation}
% where $\hat{F}_{n,k}(t) = n^{-1}\sum_{j=1}^n \1_{\{X_{j,k}\leq t\}}$ denotes the empirical cumulative
% distribution function (c.d.f.) of the $k$th coordinate
% $X_{1,k}$. The matrix $\hat\Sigma_n$ thus corresponds to an empirical 
% estimate of the covariance matrix $\Sigma$ with general term
% \begin{equation}\label{eq:def_Sigma}
% \Sigma_{k k'} = \cov\left(F_{k}(X_{1,k});F_{k'}(X_{1,k'})\right),\; 1\leq k,k'\leq K \;,
% \end{equation}
% $F_k$ denoting the c.d.f.\ of $X_{1,k}$ which we shall assume in
% the sequel to be continuous.
Theorem~\ref{th:MKW}, proved in Appendix~\ref{sec:proofMKW}, describes the
limiting behaviour of the test statistic $T(n_1,\dots,n_{L-1})$ under the null hypothesis.

%\fxnote{Would like to replace $\lambda$ by $t$ for consistency}
\begin{theorem}
\label{th:MKW}
Assume that $(\mathbf{X}_i)_{1\leq i\leq n}$ are $\mathbb{R}^{K}$-valued
  i.i.d.\ random vectors such that, for all $k$, the c.d.f.\ $F_{k}$
  of $X_{1,k}$ is a continuous function. Assume also that for 
  $\ell=0,\dots,L-1$, there exists $t_{\ell+1}$ in $(0,1)$ such
  that $(n_{\ell+1}-n_\ell)/n\to t_{\ell+1}$, as $n$ tends
  to infinity.
Then, $T(n_1,\dots,n_{L-1})$ defined in (\ref{eq:KW_mult}) satisfies
 \begin{equation}
    \label{eq:KW_mult_limite}
    T(n_1,\dots,n_{L-1})\inlaw  \chi^2\left((L-1)K\right)\;,\textrm{ as } n\to\infty\;,
  \end{equation}
 where $d$ denotes convergence in distribution and
$\chi^2((L-1)K)$ is the chi-square distribution with $(L-1)K$
degrees of freedom.
\end{theorem}

Observe that \eqref{eq:KW_mult} extends 
the classical Kruskal-Wallis test used for univariate observations to the multivariate setting. 
Indeed, when $K=1$, \eqref{eq:KW_mult} is equivalent to
\begin{equation}\label{eq:KW}
T(n_1,\dots,n_{L-1})=\frac{12}{n^2}\sum_{\ell=0}^{L-1}(n_{\ell+1}-n_\ell)\left(\bar{R}_{\ell}^{(1)}-(n+1)/2\right)^2\;,
\end{equation}
where we have replaced $\hat{\Sigma}_{n,11}$ by $\Sigma_{11} = 4\var(F_1(X_{1,1})) = 4\var(U)=1/3$ ($U$ denoting a uniform random variable on
$[0,1]$).
In the case where there is only one change-point, \textit{i.e.}, when
$L=2$, \eqref{eq:KW_mult} reduces to the test statistic proposed in 
Section~\ref{sec:test_twosample}. Indeed, using (\ref{eq:ranks}),
$T(n_1)$ can be rewritten as follows
\begin{multline*}
T(n_1)=\frac{n n_1 (n-n_1)}{n^2
  n_1}\mathbf{U}_{n}(n_1)'\hat{\Sigma}_{n}\mathbf{U}_{n}(n_1)
+\frac{n n_1
  (n-n_1)}{n^2(n-n_1)}\mathbf{U}_{n}(n_1)'\hat{\Sigma}_{n}\mathbf{U}_{n}(n_1)\\
=\mathbf{U}_{n}(n_1)'\hat{\Sigma}_{n}\mathbf{U}_{n}(n_1)=S_n(n_1)\;,
\end{multline*}
where $S_n(n_1)$ is defined in (\ref{eq:StatWL}).

\subsection{Power of the homogeneity test in the two sample case}
\label{sec:power}
In this section, we focus on the two sample homogeneity test statistic
$S_n(n_1)$ defined in~(\ref{eq:StatWL}) using local alternatives consisting of
multivariate shifts to investigate the statistical power of the proposed
approach.  We consider the following alternative hypotheses $(H_{1,n})$:
``$(\mathbf{X}_{1}, \dots, \mathbf{X}_{n_1})$ are i.i.d with common multivariate
p.d.f. $f(\mathbf{x})$ and $(\mathbf{X}_{n_1+1}, \dots, \mathbf{X}_{n})$ are
i.i.d with density $f(\mathbf{x}-\boldsymbol{\delta}/\sqrt{n})$'', 
where $f$ is symmetric, positive, and, continuously differentiable and $\boldsymbol{\delta}$ denotes an arbitrary
shift vector in $\mathbb{R}^k$. Our results can thus be compared directly with those
obtained by %\cite{mottonen:oja:tienari:1997},
\cite{hettmansperger:mottonen:oja:1998}, \cite{oja:1999} and
\cite{topchii:tyurin:oja:2003} for affine invariant multivariate generalisations of rank
tests.

We first recall the classical result pertaining to the Hotelling-$T^2$ test
\[
  H_n(n_1)= (n_1(n-n_1)/n)(\mathbf{\bar{x}_{n_1}}-\mathbf{\bar{x}_{n-n_1}})'\hat{C}_n^{-1}
(\mathbf{\bar{x}_{n_1}}-\mathbf{\bar{x}_{n-n_1}}) \; ,
\]
where $\mathbf{\bar{x}_{n_1}}=(\sum_{i=1}^{n_1}\mathbf{X}_i)/n_1$,
$\mathbf{\bar{x}_{n-n_1}}=(\sum_{i=n_1+1}^{n}\mathbf{X}_i)/(n-n_1)$,
and, $\hat{C}_n$ is the empirical covariance matrix of the
$\mathbf{X}_{i}$'s. Under $(H_{1,n})$, assuming that $n$ tends to infinity with $n_{1}/n \to
 t_1$ and under appropriate moment conditions, it holds that $H_n(n_1) \stackrel{d}{\longrightarrow} \chi_K^2(d_H(\boldsymbol{\delta}))$, where
 \begin{equation}
   \label{eq:H-T2:local_alt}
   d_H(\boldsymbol{\delta}) = t_1(1-t_1)\boldsymbol{\delta}'C^{-1}\boldsymbol{\delta} \; ,   
 \end{equation}
 $\chi_K^2(d)$ denoting the non-central chi-squared distribution with
 $K$ degrees of freedom and non-centrality parameter $d$
and $C$ denoting the covariance matrix of the
$\mathbf{X}_{i}$'s, see \cite{bickel:1965}.

Let $F(\mathbf{X_i})=(F_1(X_{i,1}),\dots,F_K(X_{i,K}))'$ denote the $K$-dimensional vector of marginal distribution functions and $\nabla\log f(\mathbf{X}_i)$ the score function. The following theorem (proved in Appendix
\ref{appendix:power}) establishes the asymptotic behaviour of $S_n(n_1)$ under $(H_{1,n})$.

\begin{theorem}\label{th:power}
Assume that
$\mathbf{X}_{1},\dotsc,\mathbf{X}_{n_{1}},\mathbf{X}_{n_{1}+1},\dotsc,\mathbf{X}_{n}$
are $\Bbb R^{K}$-valued random vectors distributed under $(H_{1,n})$ and that the $K\times K$ covariance matrix $\Sigma$
 defined by \eqref{eq:CovMatrix}
is positive definite. Assume also that the Fisher information matrix
$I_f=\PE_f[\nabla\log f(\mathbf{X}_1)\nabla\log f(\mathbf{X}_1)']$
is finite and that the densities $f_k$ of $X_{1,k}$ are upper bounded for all $k$.
Then, as $n$ tends to infinity with $n_{1}/n \to
 t_1 \in (0,1)$,
\begin{equation}\label{eq:conv_U}
%\mathbf{U}_n(n_1)'\Sigma^{-1}\mathbf{U}_n(n_1)
S_n(n_1)
\stackrel{d}{\longrightarrow}
\chi_K^2\left(4t_1(1-t_1)\boldsymbol{\delta}'A'\Sigma^{-1}A\boldsymbol{\delta}\right)=\chi_K^2(d_S(\boldsymbol{\delta}))\;,
\end{equation}
where $\Sigma$ is defined in~(\ref{eq:CovMatrix}) and $A=\PE_f[(F(\mathbf{X}_1)-1/2)\nabla\log f(\mathbf{X}_1)']$.
\end{theorem}

The following corollary, which is proved in Appendix \ref{appendix:coro:indep}, 
particularises this results to the case where the coordinates of the observations are independent.

% Always true (?) \textcolor{red}{ and that $\lim_{x\to\infty} f(x)=0$}
\begin{coro}
\label{coro:indep}
Assume that the density $f$ may be written as the product of its
marginals, $f(\mathbf{x})=\prod_{k=1}^K f_k(x_k)$. Then,
$$
d_H(\boldsymbol{\delta}) = t_1(1-t_1)\sum_{k=1}^K\frac{\delta_k^2}{\sigma_k^2}
\; \;\text{and} \; \; 
d_S(\boldsymbol{\delta})=12 t_1(1-t_1)\sum_{k=1}^K\frac{\delta_k^2}{\sigma_k^2}
\lambda_k^2\;,
$$
where $\sigma_k^2= \int x^2 f_k(x)\rmd x$ and $\lambda_k=\int \left(\sigma_k f_k(\sigma_k x)\right)^2\rmd x$.

Moreover, $d_S(\boldsymbol{\delta}) \geq (108/125) d_H(\boldsymbol{\delta})$, with $d_S(\boldsymbol{\delta})$ being equal to $(3/\pi) d_H(\boldsymbol{\delta})$ when $f_k$ are Gaussian densities.
\end{coro}

Corollary~\ref{coro:indep} states that the so-called {\emph asymptotic
  relative efficiency (ARE)} with respect to the Hotelling-$T^2$ test
$d_S(\boldsymbol{\delta})/d_H(\boldsymbol{\delta})$ is lower bounded by
$108/125\approx 0.86$. Granted that this result holds irrespectively of the
choice of the marginals $f_k$ this is a strong guarantee that extends the
well-known scalar result. As will be illustrated in Section~\ref{sec:outlier},
the ARE is much larger whenever some of the marginal have strong tails. In the
particular case of Gaussian densities, the ARE of the proposed test is equal to
$3/\pi\approx 0.95$, which is comparable to the values obtained for affine
invariant multivariate rank tests \cite[p. 331]{oja:1999}.

For the particular case of multivariate Gaussian distributions, an explicit expression of the ARE
can be obtained without assuming independence, as shown by the following corollary (proved in
Appendix \ref{appendix:coro:gaussian}).

\begin{coro}
\label{coro:gaussian}
Assume that $f$ is a multivariate Gaussian p.d.f. with mean 0 and covariance matrix $C$, then the matrices $A$ and $\Sigma$ featured in Theorem ~\ref{th:power} are given by $A = (2\sqrt{\pi})^{-1}\operatorname{diag}(\sigma_1^{-1},\dots, \sigma_K^{-1})$ and $\Sigma_{k\ell} = 2/\pi \arcsin(C_{k\ell}/(2\sigma_k\sigma_\ell))$, for $1 \leq k, \ell \leq K$, where $(\sigma_k^2)_{1\leq k \leq K}$ denote the diagonal elements of $C$.

Furthermore, the asymptotic relative efficiency of the test statistic $S_n(n_1)$ can be lower bounded as follows
$$
d_S/(\boldsymbol{\delta})d_H(\boldsymbol{\delta})\geq (3/\pi) (\sigma_{\min}^{2}/\sigma_{\max}^{2})(\lambda_{\min}(C)/\lambda_{\max}(|C|)) \; ,
$$
where $\sigma_{\min}^2$ and $\sigma_{\max}^2$ denote, respectively, the minimal
and maximal diagonal terms of $C$, $\lambda_{min}(C)$ is the minimal eigenvalue
of $C$ and $\lambda_{max}(|C|)$ the maximal eigenvalue of
$|C|=(|C_{k,\ell}|)_{1\leq k,\ell\leq K}$.
\end{coro}

We observe here an important difference with the test statistic of
\cite{mottonen:oja:tienari:1997} which is designed so as to guarantee that its
ARE in the multivariate Gaussian case does not depend on the value of $C$. For
the proposed statistic, the ARE is lower bounded by the minimal eigenvalue of
$A'\Sigma^{-1}AC$, where $A$ and $\Sigma$ are defined in
Corollary~\ref{coro:gaussian}. Recall that $\Sigma$ is proportional to the
Spearman correlation matrix whereas Corollary~\ref{coro:gaussian} implies that
$A C A'$ is proportional to the standard correlation matrix. Empirically, it
can be checked using numerical simulation that the minimal eigenvalue of
$A'\Sigma^{-1}AC$ can only be small when $C$ itself is poorly conditioned. The
second statement of Corollary~\ref{coro:gaussian} substantiates this claim by
providing a lower bound which, for positively correlated $C$ at least, is
inversely proportional to the condition number of $C$. Note that this bound
appears to be pessimistic in practice as for values of $K$ in the range 2 to 6,
we observed the ARE to be larger than 0.8 for all matrices with condition
number smaller than 10 (recall that the ARE is equal to 0.95 in the scalar
case); for matrices with condition number up to 100, the ARE was still larger
than 0.3. Hence, in situations where the coordinates of the data are not too
correlated, despite the fact that the proposed test is not affine invariant,
its loss with respect to the (optimal) likelihood ratio test in the Gaussian
case is usually negligible, a fact that will be illustrated by the numerical
simulations of Section~\ref{sec:Eval1}.

\section{Change-point estimation and detection}
\label{sec:ChangePoint}

We now consider the setting in which the position of the
potential change-points are unknown (still assuming that their number is known). Our proposal is to consider the statistic
$T(n_1,\dots,n_{L-1})$ described in the previous section and to optimise it
over all the possible change point locations. This proposal is however
faced with two serious difficulties. The first one, which is of
computational nature is related to the feasibility of the maximisation
when there is more than a single change-point. We start by showing
that the maximisation of $T(n_1,\dots,n_{L-1})$ is amenable to dynamic
programming and stays feasible even when $L$ is large. The second
difficulty, to which a partial answer is provided in
Section~\ref{sec:CP-onechange}, is statistical and concerns the
interpretation of the value of $T(n_1,\dots,n_{L-1})$ as optimising
with respect to the change-point location obviously modifies the
distribution of the values of the test statistic. This is a difficult
issue in general, but we show how to obtain meaningful and simple-to-compute $p$-values for a variant of the test in the case of a single change-point.

\subsection{Multiple change-point estimation}
\label{sec:CP-multiple}
Assuming a known number of change-points $L$, we propose to use the test statistic
described in Section~\ref{sec:test_homogeneity_multigroups} to determine the positions of the segment boundaries
$n_1,\dots,n_{L-1}$. These unknown change-point locations are estimated by maximising the
statistic $T(n_1,\dots,n_{L-1})$ defined in (\ref{eq:KW_mult}) with respect to $n_1,\dots,n_{L-1}$:
\begin{equation}\label{eq:est_chang}
(\hat{n}_1,\dots,\hat{n}_{L-1})=
\argmax\limits_{1\leq n_1<\dots<n_{L-1}\leq n} T(n_1,\dots,n_{L-1})\;.
\end{equation}
% We shall use $\hat{t}_1,\dots,\hat{t}_{K^{\star}-1}$ as change-point
% estimators as soon as the hypothesis $(H_0)$ is rejected.

In practice, direct maximisation by enumeration in (\ref{eq:est_chang}) is computationally
prohibitive as it corresponds to a combinatorial task whose complexity grows exponentially with
$L$. However, due to the fact that the matrix $\hat{\Sigma}_n$ is common to all segments, the
statistic $T(n_1,\dots,n_{L-1})$ defined in (\ref{eq:KW_mult}) has an additive structure which
makes it possible to adopt a dynamic programming strategy. We refer here to the classical dynamic
programming approach to the segmentation task which is described in \cite{kay:1993} used by, among
others, \cite{bai:perron:2003} and can be traced back to the note by \cite{Bellman:1961}. More
precisely, using the notations
$$
\Delta(n_{\ell}+1:n_{\ell+1})=(n_{\ell+1}-n_\ell)\mathbf{\bar{R}}_\ell'\;\hat{\Sigma}_n^{-1}\;
\mathbf{\bar{R}}_\ell\;,
$$
and
$$
I_{L}(p)=\max_{1<n_1<\dots<n_{L-1}<n_{L}=p}\sum_{\ell=0}^{L-1}\Delta(n_{\ell}+1:n_{\ell+1})\;,
$$
we have
\begin{equation}\label{eq:dyn_recurs}
I_{L}(p)=\max_{n_{L-1}}\left\{I_{L-1}(n_{L-1})+\Delta(n_{L-1}+1:p)\right\}\;.
\end{equation}
Thus, for solving the optimisation problem (\ref{eq:est_chang}), we
proceed as follows. We start by computing the $\Delta(i:j)$ for all
$(i,j)$ such that
$1\leq i<j\leq n$. All the $I_1(E)$ are thus
available for $E=2,\dots,n$. Then $I_2(E)$ is computed by using
the recursion (\ref{eq:dyn_recurs}) and so on. The overall numerical
complexity of the procedure is thus proportional to $L \times n^2$ only.

\subsection{Assessing the significance of the test in the single change-point case}
\label{sec:CP-onechange}

In addition to practical algorithms for estimating change-point locations, one needs tools to
assess the plausibility of the obtained change-point configuration. An important step in that
direction is to characterise the behaviour of $T(\hat{n}_1,\dots,\hat{n}_{L-1})$ under the null
hypothesis that the data are indeed fully homogeneous. This is a difficult issue in general due to
the optimisation over all possible change-point configurations. A possible calibration approach
consists in running Monte Carlo experiments, possibly using bootstrap techniques if a representative
sample of the baseline data of interest is available. We show below that in the case where $L=2$,
\textit{i.e.}, when looking for a single potential change-point, it is possible to obtain a simple
computable approximation to the asymptotic $p$-value of the test.

To do so, we consider in the rest of this section a modification of the test statistic used in
(\ref{eq:est_chang}). The practical consequences of using this variant rather than the statistic
$T(\hat{n}_1)$ when $L=2$ will be discussed after Theorem~\ref{theo:change-point}
which states the main result of this section.

Let $\mathbf{V}_{n}(n_1)=(V_{n,1}(n_1),\dots, V_{n,K}(n_1))'$ denote the vector such that 
\begin{equation}\label{eq:def:V_n}
V_{n,k}(n_1)= \inv{n^{3/2}} \sum_{i=1}^{n_1} \sum_{j=n_1+1}^{n} \left\{
   \1(X_{i,k}\le X_{j,k})-\1(X_{j,k}\le X_{i,k}) \right\},\; k=1,\dotsc,K\; ,
\end{equation}
and define
\begin{equation}\label{eq:S_n_tilde}
\tilde S_{n}(n_1) = \mathbf{V}_{n}(n_1)' \hat{\Sigma}_n^{-1} \mathbf{V}_{n}(n_1)\;.
\end{equation}
%as soon as the $K\times K$ matrix $\hat{\Sigma}_n$ defined in
%(\ref{eq:Sigma_emp}) is invertible.
Note that $\mathbf{V}_{n}$ only differs from $\mathbf{U}_{n}$ by the
normalisation, which is now independent of $n_1$. We now consider the statistic
\begin{equation}
  \label{eq:statRuptureWL}
  W_n=\max_{1\leq n_1 \leq n-1} \tilde S_{n}(n_1)\;.
%W_n=\max_{0< s<1} |\mathbf{U}_{n}(\lfloor t\cdot n \rfloor) \cdot\Gamma^{-1}\cdot\mathbf{U}_{n}(\lfloor t\cdot n \rfloor)'|\; .
\end{equation}

The following theorem, proved in Appendix~\ref{sec:Appendix2}, gives the asymptotic $p$-values of
$W_n$ under the null hypothesis that no change in distribution occurs within the observation data.

\begin{theorem}\label{theo:change-point}
  Assume that $(\mathbf{X}_i)_{1\leq i\leq n}$ are $\Bbb R^{K}$-valued
  i.i.d.\ random vectors such that, for all $k$, the c.d.f.\ $F_{k}$
  of $X_{1,k}$ is a continuous function. Further assume that the
  $K\times K$ matrix $\Sigma$ defined in (\ref{eq:CovMatrix}) is
  invertible. Then,
 \begin{equation}
    \label{eq:WLrupture}
    W_{n}\inlaw  \sup_{0<t<1} \left(\sum_{k=1}^{K} B_{k}^{2}(t)\right),\textrm{ as } n\to\infty\;,
  \end{equation}
 where $d$ denotes convergence in distribution and
$\{B_{k}(t),\ t\in(0,1)\}_{1\leq k\leq K}$ are independent Brownian bridges.
\end{theorem}

To determine the $p$-value $P_{\text{val}}(W_{n})$ associated to
(\ref{eq:WLrupture}), one can use the following 
result due to \cite{kiefer:1959}:
\begin{multline}
  \label{eq:pvalKiefer}
  P_{\text{val}}(b) = \P\left(\sup_{0<t<1}\left(\sum_{k=1}^{K} B_{k}^{2}(t)\right)>b\right)\\
=1-  \frac{4}{\Gamma(\frac{K}{2})2^{\frac{K}{2}}b^{\frac{K}{2}}}
  \sum_{m=1}^{\infty}
  \frac
  {(\gamma_{(K-2)/2,m})^{K-2}\exp[-(\gamma_{(K-2)/2,m})^2]/2b}
  {[J_{K/2}(\gamma_{(K-2)/2,m})]^{2}}\;,
\end{multline}
where $J_{\nu}$ is the Bessel function of the first kind,
$\gamma_{\nu,m}$ is the $m$-th nonnegative zero of $J_{\nu}$ and
$\Gamma$ is the Gamma function. 
In practice, only a few terms of the series have to be computed.
% ,
% although higher values of the dimensions $K$ and/or of the
% statistic $b$ require a larger number of terms
For values of $K$ of forty or less computing the $p$-values from the thirty first terms of the series was sufficient.

As noted at the beginning of Section~\ref{sec:CP-onechange}, the normalisation of $\mathbf{V}_{n}$
differs from that of $\mathbf{U}_{n}$, resulting in a statistic $W_n$ that does not coincide with
$T(\hat{n}_1)$. From our practical experience, replacing $\mathbf{V}_{n}$ by $\mathbf{U}_{n}$ in
the definition of $W_n$, that is using $T(\hat{n}_1)$ instead of $W_n$, produces a statistic that
has the same detection and localisation capacities when the potential change occurs in the central
region of the observation window, say between $n/4$ and $3n/4$ observations. For potential changes
occurring closer to the beginning or to the end of the observation window, $T(\hat{n}_1)$ has an enhanced
detection power at the expense of a slight increase in the rate of false alarms, with corresponding
spurious detections occurring mostly near the borders of the observation window. Proceeding as in
Appendix~\ref{sec:Appendix2}, one can prove a result related to Theorem~\ref{theo:change-point} for
$T(\hat{n}_1)$ (used when $L=2$) by imposing some additional conditions on the admissible values
of $n_1$ (namely, that the maximum is searched only for value of $n_1$ such that $n_1/n$ is bounded
from above and below). The resulting limit, expressed in terms of Bessel
processes, does not yield easily computable asymptotic $p$-values, although approximations such
as those studied by \cite{estrella:2003} could be used for approximating extreme quantiles (that
is, very low $p$-values).

\section{Numerical Experiments}
\label{sec:Eval1}

In this section, we report numerical experiments that illustrate different aspects of the
methods proposed in Sections~\ref{sec:TestingHomo}
and~\ref{sec:ChangePoint}. 
A software implementing these methods in Python is available as a
supplementary material of the paper.
For easy reference,
the two- or multi-sample homogeneity test defined by~(\ref{eq:KW_mult}) is referred to as
 \emph{MultiRank-H} in the following; the change point estimation criterion defined in~(\ref{eq:est_chang}) is 
referred to as \emph{MultiRank}.

\subsection{Illustration of the two-sample homogeneity test}
\label{sec:simulToy}

\begin{figure}[hbt] \centering %multirank/code/scripts/star/resStar.py
\begin{tabular}{cc}
  \multicolumn{2}{c}{\includegraphics[width=.48\textwidth]{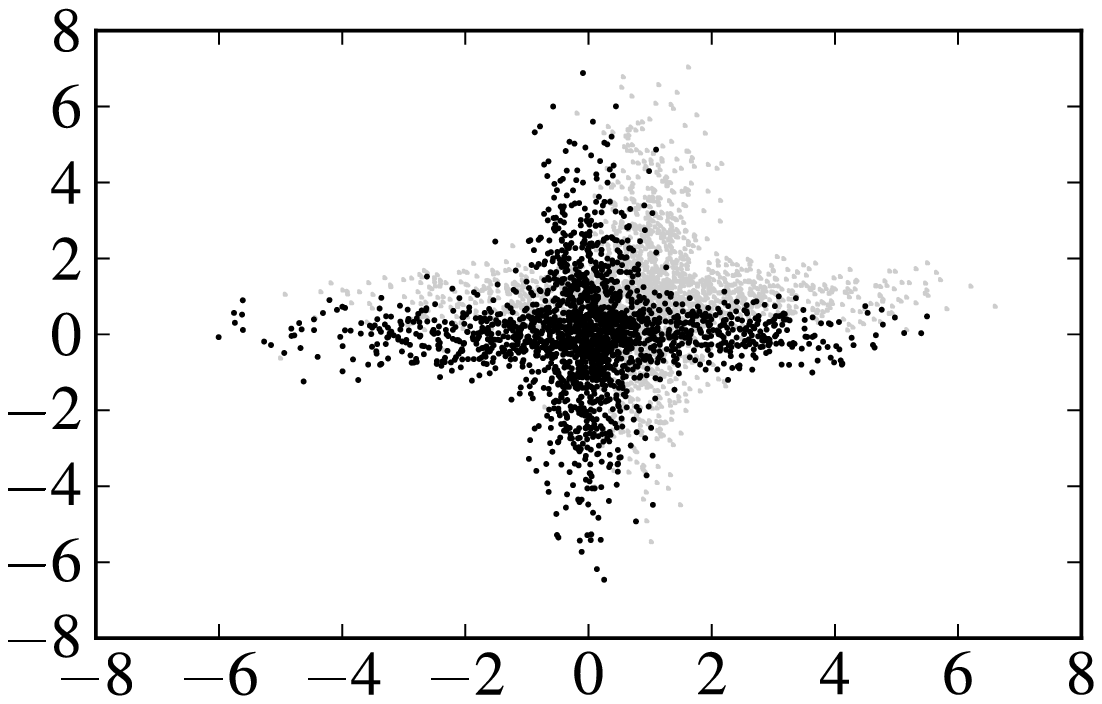}}\\
  \multicolumn{2}{c}{(a)}\\
  \includegraphics[width=.48\textwidth]{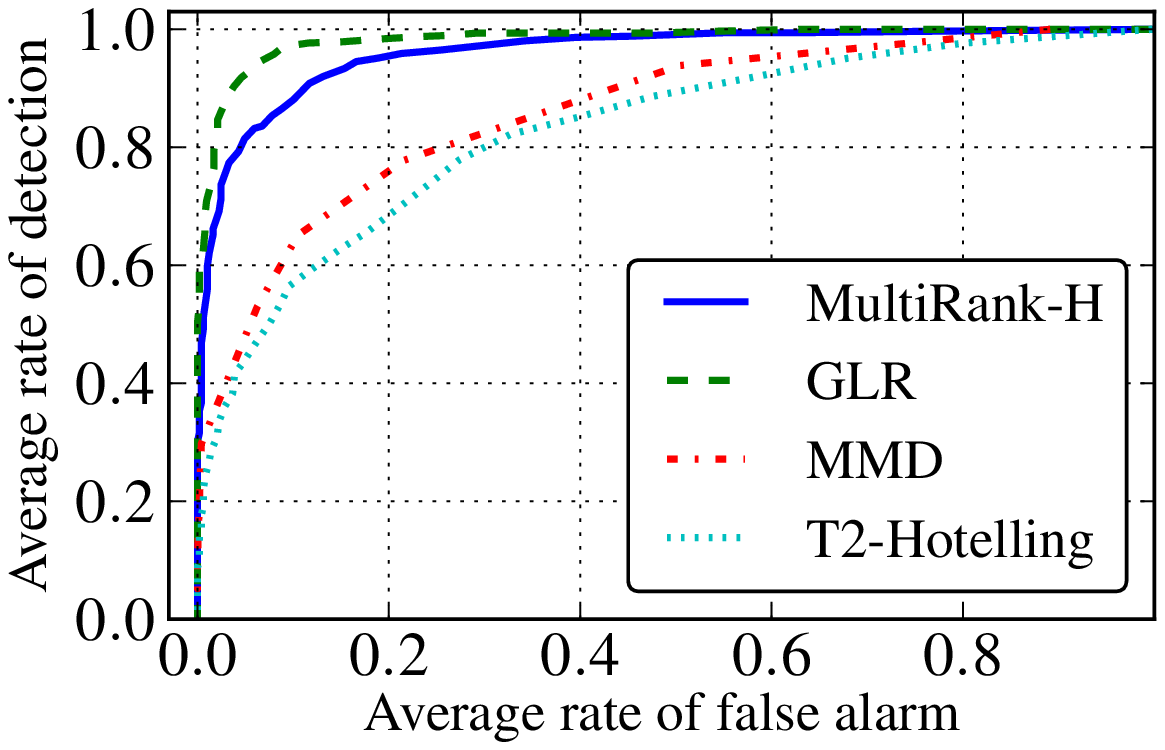} &
  \includegraphics[width=.48\textwidth]{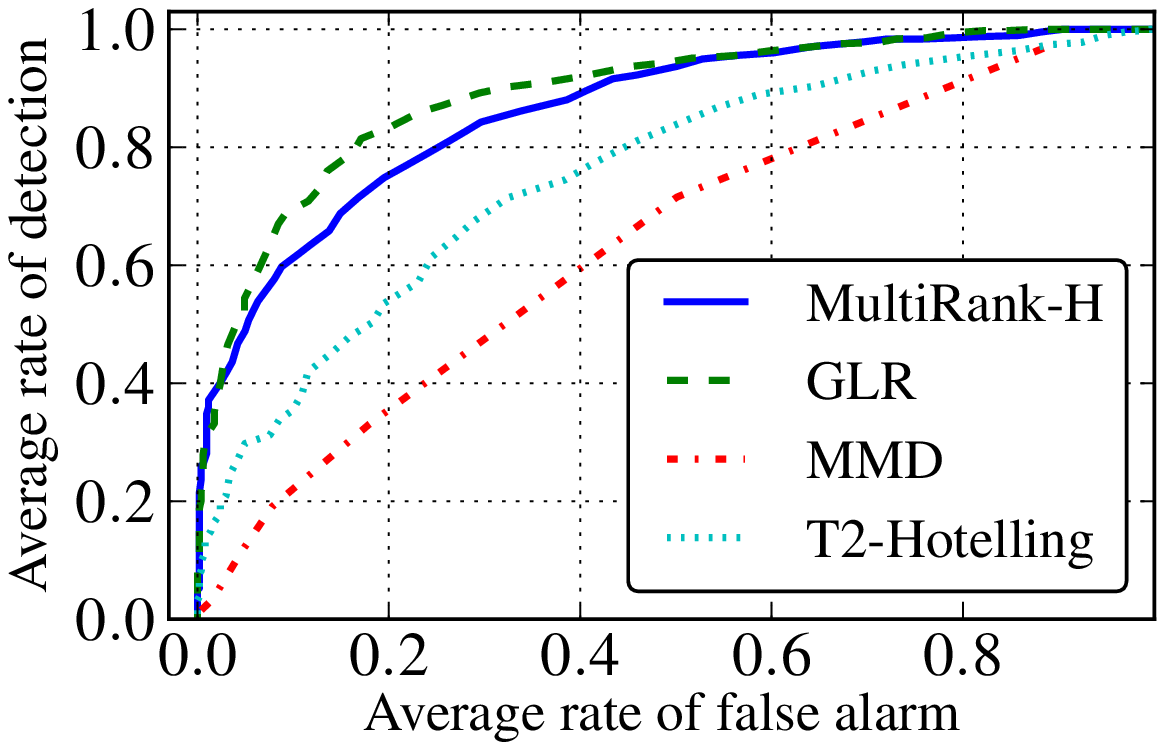} \\
  (b) & (c) 
\end{tabular}
  \caption{(a) Example of observations under baseline and
alternative distributions, (b) ROC curves for MultiRank-H, MMD, Hotelling's $T^2$ and
LR, (c) Same as (b) with eight-dimensional Gaussian noise padding (that is, with $K=10$).}
  \label{fig:Cross}
\end{figure}

We start by considering the basic two-sample
homogeneity test first introduced in Section~\ref{sec:test_twosample}.
For this, we generate baseline observations distributed
as a mixture of two two-dimensional Gaussian
densities with common mean $(0,0)$ and diagonal covariance matrices
with diagonal terms equal to $(4, 0.2)$ and $(0.2, 4)$, respectively.
For the alternative distribution, we generate observations having the same
characteristics except that the mean is now equal to $\boldsymbol{s}=(0.5,0.5)$. In this case,
$n=100$ and the two groups (baseline and alternative data) are of length $n/2=50$.
%\fixme{Please check}
Figure~\ref{fig:Cross} (a) shows a typical example of the data, represented
as a two-dimensional scatter plot.

MultiRank-H is compared with three other approaches. The first is the
the Maximum Mean Discrepancy (MMD) 
statistics proposed by \cite{gretton:2006:nips}, which is a kernel-based
test here used with a Gaussian kernel having a
bandwidth given by the median distance between the samples as suggested by
\cite{gretton:2006:nips}, \cite{desobry:davy:doncarli:2005} and \cite{harchaoui:bach:moulines:2008}.
The second approach is the classical Hotelling's $T^2$-test \citep[p. 67]{chen:gupta:2000} which
is optimal in the multivariate Gaussian case. The third method is to use the
likelihood ratio (LR) test assuming a known structure for the model whose
parameters (mean vectors and diagonal terms of the covariance matrices) are
estimated using the Expectation-Maximisation algorithm.
This latter approach is optimal in this context but is the only one that 
uses some knowledge about the distribution of the data.
These methods are compared through their ROC (Receiver Operating
Characteristic) curves, averaged over $1000$ Monte Carlo
replications of the data.

From the results of Section~\ref{sec:power}, we know that, as the covariance matrix of the data is
proportional to the identity matrix, the asymptotic performance of Hotelling's $T^2$-test does not
depend on the particular choice of the shift $\boldsymbol{s}$, but only on its $L^2$ norm. For
the Multirank-H approach the answer is less straightforward as the coordinates of the data are not
independent in this example. However, Monte Carlo estimation of the asymptotic performance index
$d_s(\boldsymbol{\delta})$ in~(\ref{eq:conv_U}) shows that it does not varies by more than 1\% with
the direction of $\boldsymbol{\delta}$. We indeed verified that in the setting of
Figure~\ref{fig:Cross} the variation in performance with the direction of the shift $\boldsymbol{s}$ is negligible.

The results displayed in Figure~\ref{fig:Cross}(b)
show that MultiRank-H is on a par with the LR
and outperforms the other two approaches. MMD performs somewhat better than 
Hotelling's $T^2$ in this context due to the non-Gaussian nature of the data.
Figure~\ref{fig:Cross}(c) corresponds to the more difficult setup in which
eight-dimensional i.i.d. Gaussian
random vectors of variance $2.5^2$ are appended to the data described previously. In this case,
the data are thus ten-dimensional but the change only affects two coordinates.
MultiRank-H is again comparable to the LR, which is still optimal in this case.
Note the lack of robustness of MMD which is distinctly dominated by Hotelling's $T^2$
in this second scenario.

% Probably not useful.
%:
% ${\cal C}_{1} \sim \frac{1}{2} {\cal N}(0,\Sigma_1) + \frac{1}{2}{\cal
%   N}(0,\Sigma_2)$ and
% ${\cal C}_{2} \sim \frac{1}{2} {\cal N}(0.5,\Sigma_1) + \frac{1}{2}{\cal
%   N}(0.5,\Sigma_2)$, with
% \begin{eqnarray*}
%   \label{eq:mix1}
%    \Sigma_1 =
%   \begin{bmatrix}
%     4&0 \\ 0&0.2
%   \end{bmatrix}
%   \text{and~}
%   \Sigma_2 =
%   \begin{bmatrix}
%     0.2 & 0 \\ 0 & 4
%   \end{bmatrix}.
% \end{eqnarray*}

\subsection{Properties of the change-point detection test}
\label{sec:jitter_ramp}

In this section, we investigate the properties of the change-point
detection test based on the statistic $T(\hat{n}_1)$
resulting from the optimisation of (\ref{eq:statRuptureWL}).
The simulations reported in this section are based on the following common
benchmark scenario: under ($H_{0}$), that is, in the absence of change, we generate 100 samples
from a five-dimensional standard Gaussian distribution. Under
($H_{1}$), the observations are similar, except that the common mean of the
Gaussian vector changes to 0.3 at a location which is either equal
to 1/4 or 1/2 of the observation window, that is, at indices 25 or 50.
This scenario corresponds to a simple situation where all coordinates of the data
possibly undergo similar upward shifts.
The ROC curves that are plotted in the following are based on 2000 replications 
of the simulated data.

\subsubsection{Comparison with marginal decisions}
\label{sec:MRvsMarg}

The MultiRank test statistic is obviously based on a combination of marginal rank
statistics. Nevertheless, it incorporates two important aspects of the multivariate change-point
detection problem: first, detection of simultaneous changes in multiple coordinates should make the
presence of an actual change-point more likely, and, second, the existence of dependence between
the coordinates should influence the decision. To illustrate these observations, we compare
MultiRank with a simpler heuristic approach that combines marginal decisions based on Bonferroni bound, using as test statistic $\max_{1\leq k\leq
  K}\max_{1\leq n_1\leq n-1} V_{n,k}(n_1)$. The results obtained with
the data-generating mechanism described at the beginning of Section~\ref{sec:jitter_ramp}
are displayed in the leftmost plot of Figure~\ref{fig:MrvsMarg}.
We also compare both approaches in a setting where the covariance
matrix of the Gaussian vector is not the identity matrix anymore but
a tridiagonal matrix with a common value of $0.45$ (positive
correlation) or $-0.45$ (negative correlation) on the sub- and super-diagonal.
The resulting ROC curves are displayed in the middle and right plots
of Figure~\ref{fig:MrvsMarg}, respectively.

The leftmost plot of Figure~\ref{fig:MrvsMarg} shows that the approach that combines the
marginal statistics by taking into account their correlation, that is
MultiRank, outperforms the Bonferroni-type approach. 
Furthermore, when the coordinates are
positively correlated, the rate of detection of the MultiRank
method decreases for a given false alarm rate and when negatively correlated, the rate of detection increases.
The performance
of the Bonferroni-type approach
on the other hand does not improve for the negatively correlated data.
The MultiRank method
captures an important feature of the problem that fails to be exploited by
the mere marginal decisions: negative correlations in the data
make the detection of simultaneous upward jumps easier while positive correlations render this task more difficult.

Although Figure~\ref{fig:MrvsMarg} deals with change-point estimation (which includes the
maximisation with respect to the change point position), we note that Corollary~\ref{coro:gaussian}
of Section~\ref{sec:power} implies that the performance of the Multirank-H homogeneity test is here
comparable to that of the Hotelling's $T^2$ test as the condition number of the covariance matrices
considered in Figures~\ref{fig:MrvsMarg} is relatively small (less than 8.1). The Multirank change detection test obviously inherits this property as its performance is nearly indiscernible from that of the LR test for all three choices of the covariance matrix in Figure~\ref{fig:MrvsMarg}.

\begin{landscape}
\begin{figure}[h]
  % multirank/code/scripts/correlated/ROC_MR_marg-res.py
  % see ROC_MR_marg_LR-res.py for figure with LR test
    \centering
    \begin{tabular}{ccc}
    \includegraphics[width=0.42\textwidth]{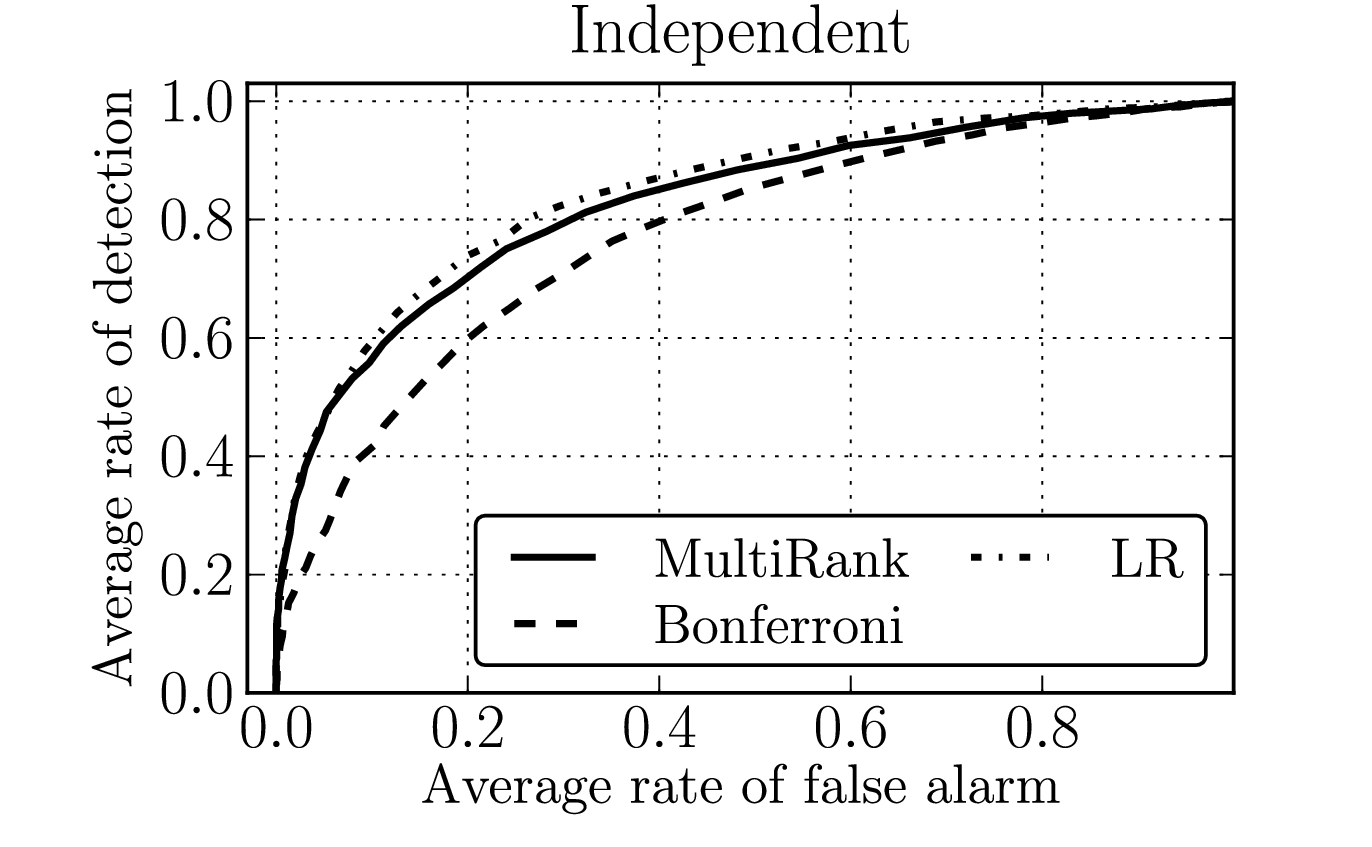}&
     \includegraphics[width=0.42\textwidth]{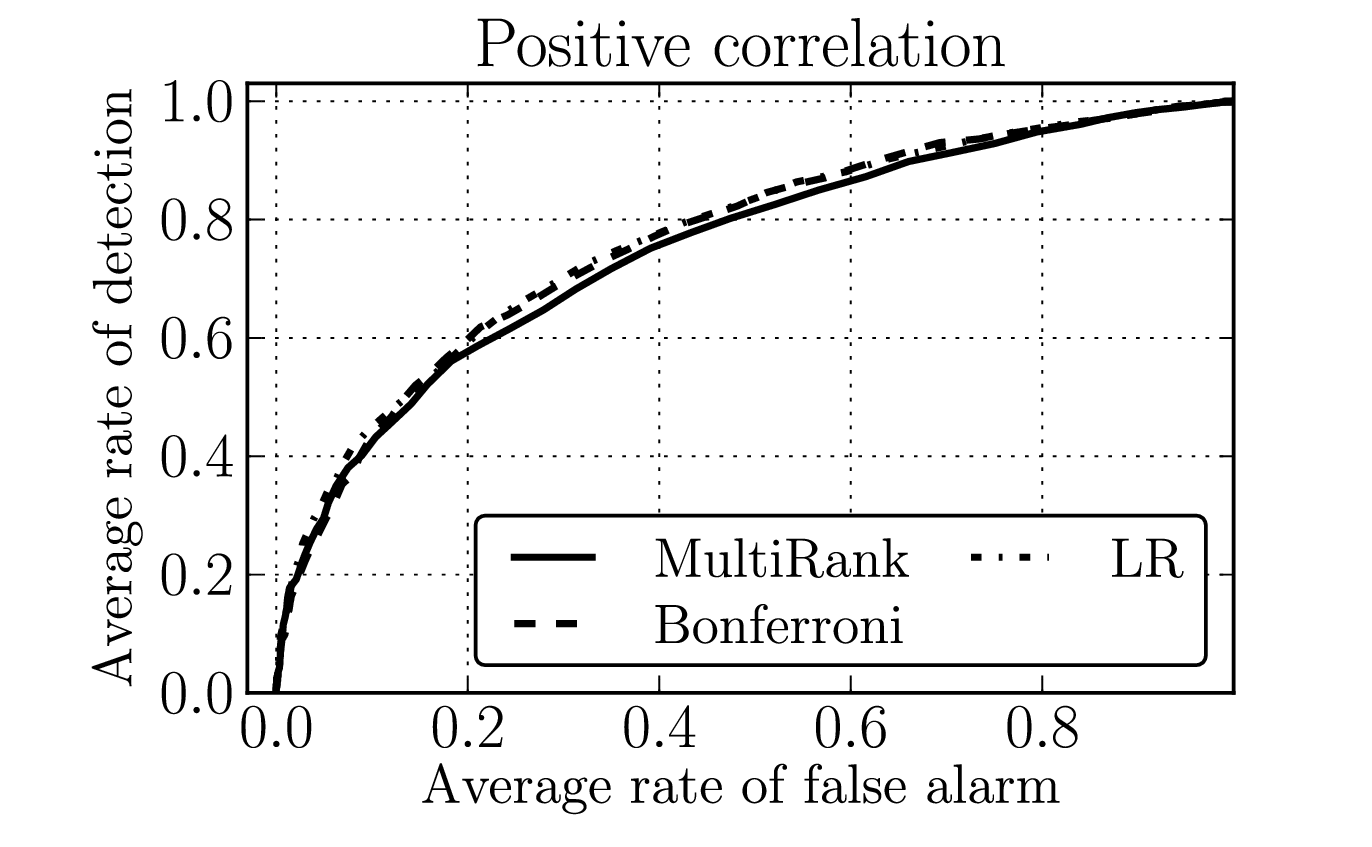}&
    \includegraphics[width=0.42\textwidth]{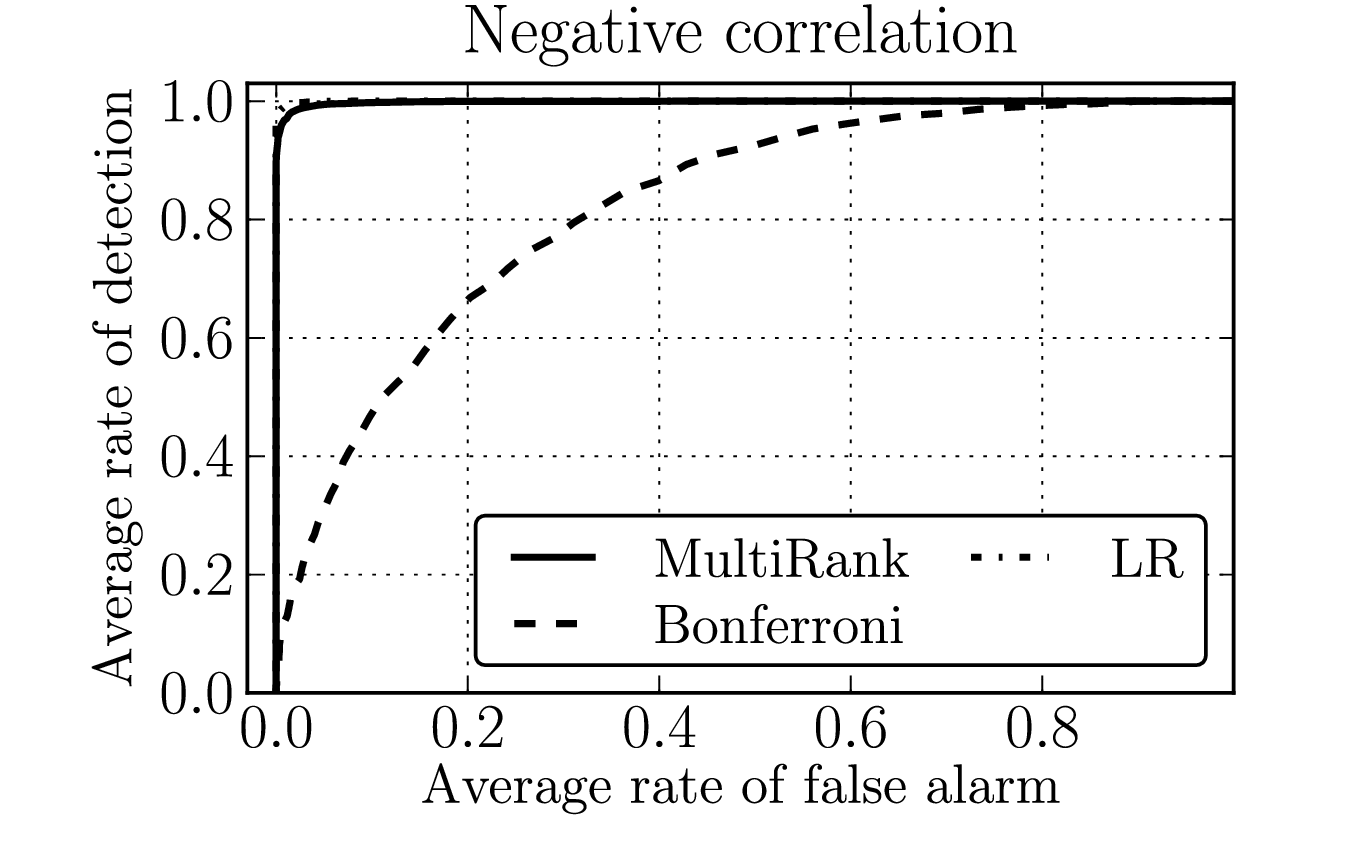}\\
  %\small{(a)}&\small{(b)}&\small{(c)}
  \end{tabular}
    \caption{ROC curves for the MultiRank and the
      Bonferroni-type approaches when the coordinates
      are independent (left), positively correlated (middle) and
      negatively correlated (right). The change-point instant is
      located at 1/4 and 1/2 of the observations window of length 500.}
    \label{fig:MrvsMarg}
\end{figure}

\begin{figure}[h]% multirank/code/scripts/outliers/besselROC_sriva_outlier-res.py
  \centering
  \begin{tabular}{ccc}
  \includegraphics[width=0.42\textwidth]{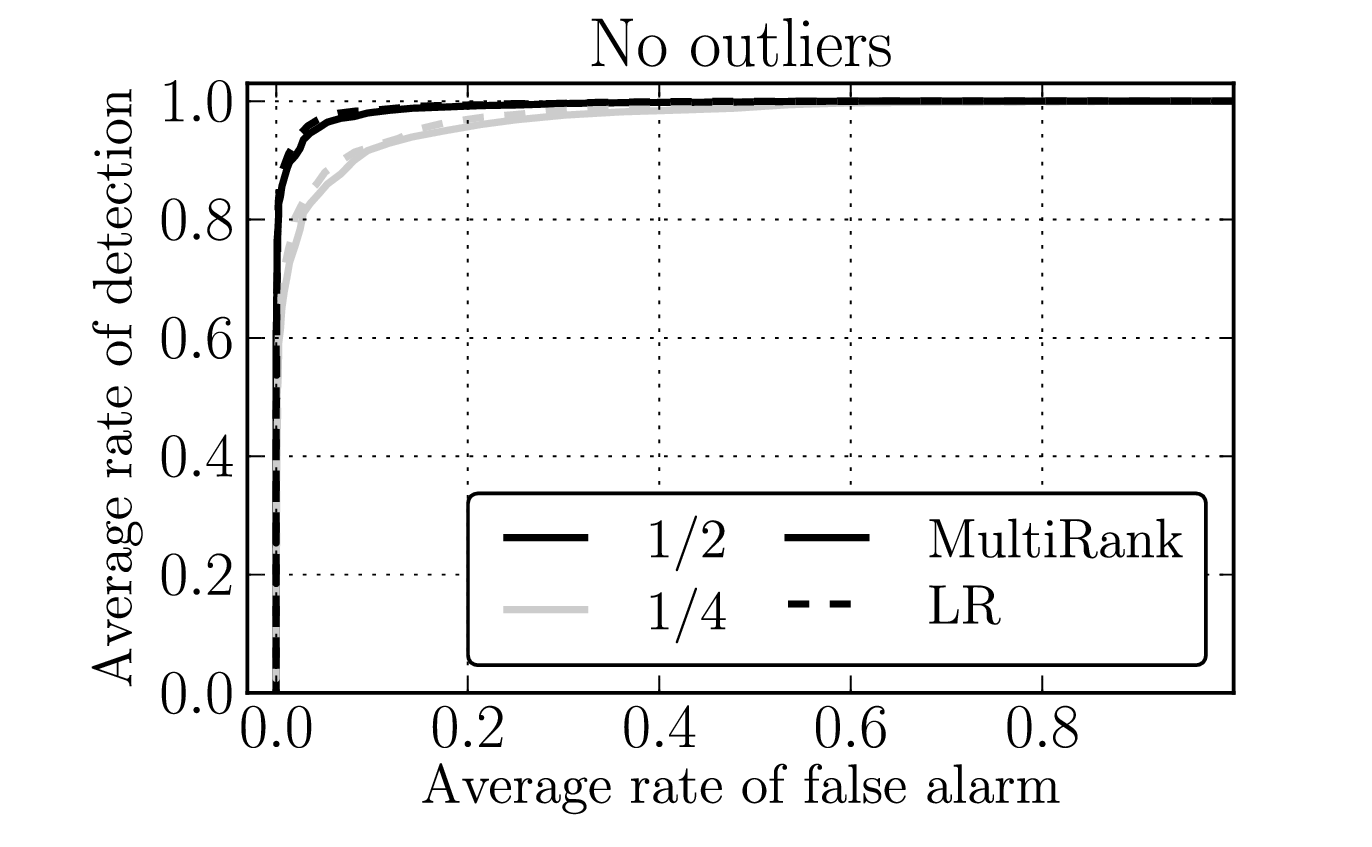}&
  \includegraphics[width=0.42\textwidth]{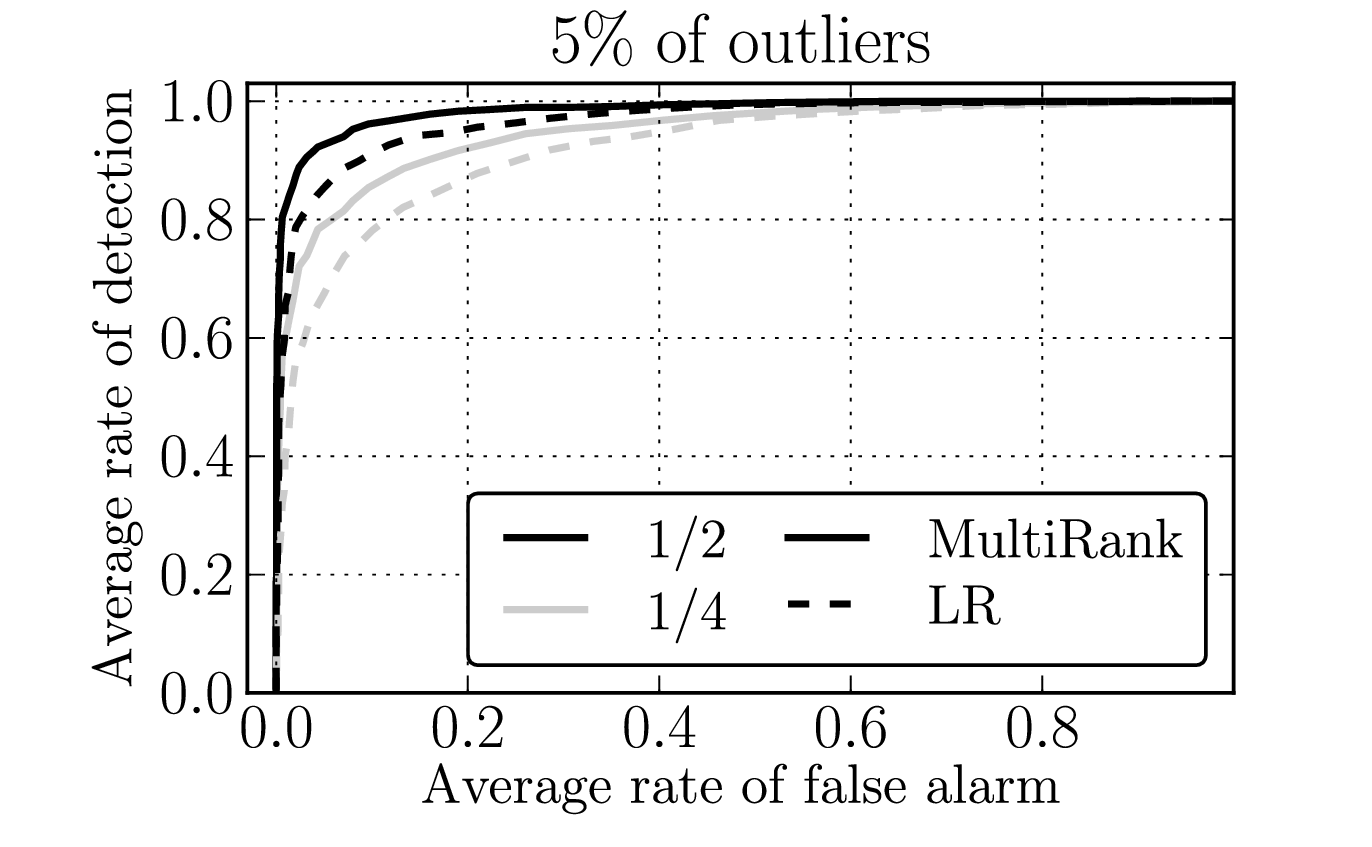}&
  \includegraphics[width=0.42\textwidth]{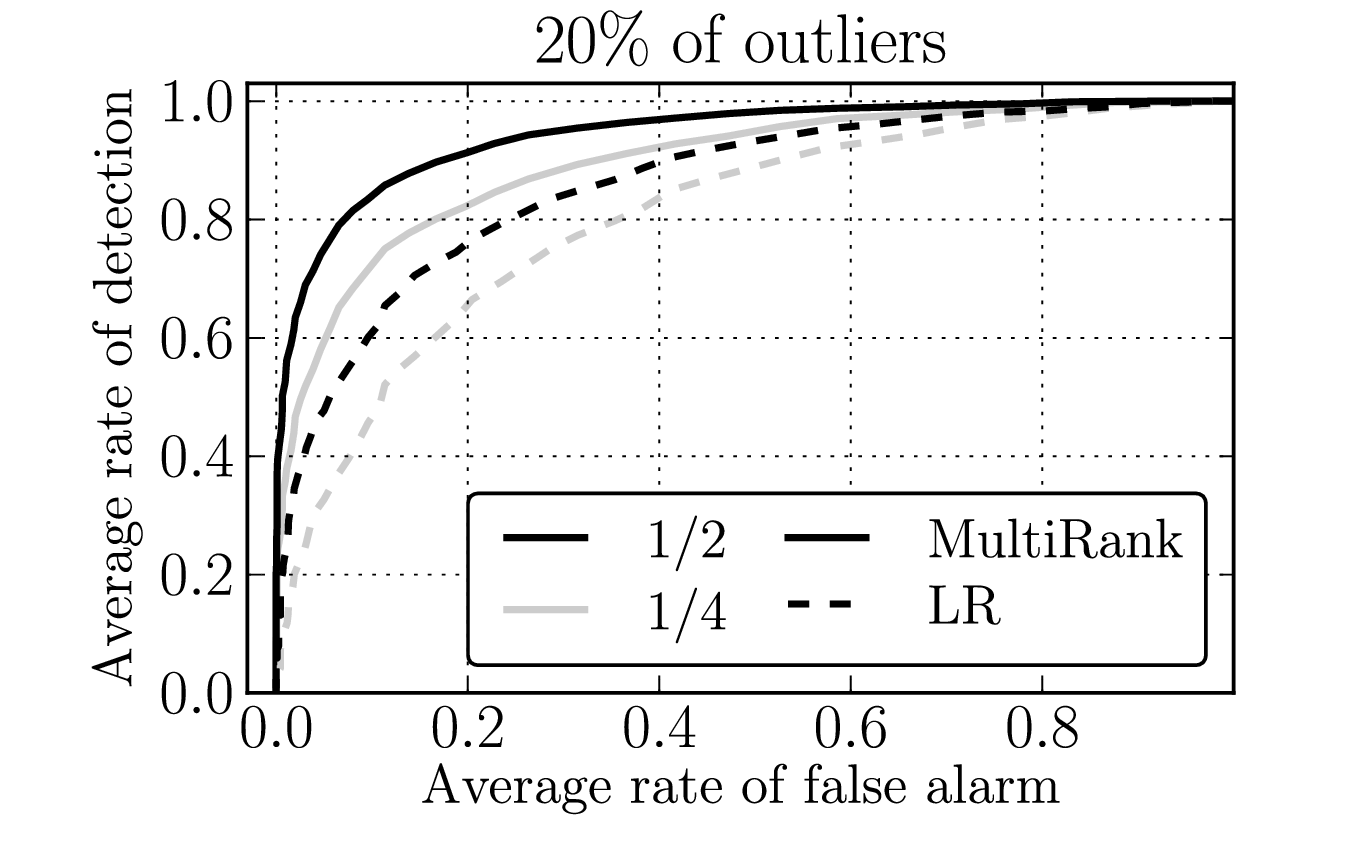}\\
  %\small{(a)}&\small{(b)}&\small{(c)}
\end{tabular}

  \caption{ROC curves for the MultiRank approach and the
    likelihood-ratio procedure (LR) for three different
    proportions of outliers 
(from left to right: 0, 5 and 20$\%$) when the change-point instant is located at 1/4 and 1/2 of the window of length 500.}
  \label{fig:ROC_outliers}
\end{figure}
\end{landscape}

\subsubsection{Robustness with respect to outliers}
\label{sec:outlier}

Here we illustrate the robustness of the MultiRank approach with respect to outliers in the data by
considering the same simulation scenario as in the previous section (with a larger shift of
amplitude 0.5) progressively contaminated by large
additive outliers. The outliers distribution is the multivariate Gaussian
distribution with covariance matrix $10\operatorname{Id}_5$, instead of
$\operatorname{Id}_5$ for the baseline distribution (where $\operatorname{Id}_5$ refers to the five by five identity matrix). The fraction of outliers is varied between 0, 5 and 20\%. The MultiRank approach is compared to the 
parametric likelihood-ratio based change-point detection test
described by \cite{srivastava:worsley:1986}. This latter method, which is itself
based on Hotelling's $T^2$-test statistic, is optimal in the absence of
outliers as the baseline and alternative distributions are both Gaussian.
As shown in the leftmost plot of Figure~\ref{fig:ROC_outliers},
MultiRank has comparable performance with the parametric approach in
the case where there are no outliers in the data. However, as shown in
the middle and rightmost plots of Figure~\ref{fig:ROC_outliers}, MultiRank
demonstrates its robustness with respect to the presence of outliers
as it barely suffers from the presence of additive outliers
contrary to the parametric approach. 

\subsection{Application to genomic hybridisation data}
\label{sec:real_data}
To illustrate the potential of the approach, we consider its application to the segmentation of multiple individual genomic data.
We consider the bladder cancer micro-array aCGH dataset %\footnote{ \url{http://cbio.ensmp.fr/~frapaport/CGHfusedSVM/index.html}} 
studied by \cite{vert:bleakley:2010} which consists of records of copy-number variations, \textit{i.e.} abnormal alteration of the quantity of DNA sections.

The objective here is to jointly segment data recorded from different subjects so 
as to robustly detect regions of frequent deletions or additions of DNA which could be characteristic of cancer.
Each of the 57 profiles provides the relative quantity of DNA for 2143
probes measured on 22 chromosomes.
We ran the change-point estimation algorithm on each of 22 chromosomes separately, thus processing 22 different 9- to 57-dimensional signals (depending on the selected groups of patients at different stages of cancer) of length 50 to 200 (the number of probes varies for each chromosome).

\begin{figure}[h]
  \centering
  %jrssbfig/rupPente.py
  \includegraphics[width=0.6\textwidth]{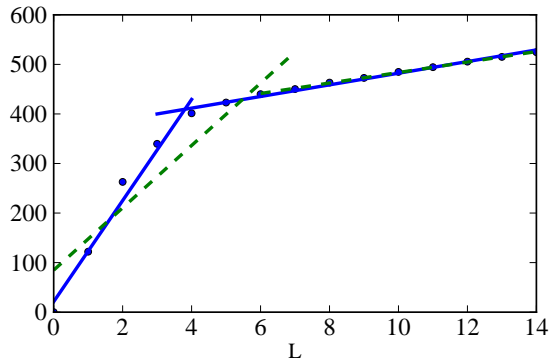}
  \caption{Determining the optimal number of change-points. Here, the actual number of change-points is $L^{\star}=4$; the optimal regression is displayed in solid lines, while a non-optimal alternative (for $L=6$) is displayed in dashed lines.}
  \label{fig:rupPente}
\end{figure}

In this paper, we have not considered principled methods for inferring the number of change-points from the data. 
We describe below an heuristic approach to determine the number of change-points which, despite its simplicity, performs in our experience much better than the use of generic penalties such as AIC or BIC.
Values of the statistics $I_L(n)$, for $L=0,\dotsc,L_{\text{max}}$, are first computed using the procedure described in Section~\ref{sec:CP-multiple}.
The algorithm is based on the principle
%, also hinted by \cite{Lavielle:2005}, 
that in the presence of $L^{\star}\geq 1$ change-points, if $I_L(n)$ is plotted
against $L$, the resulting graph can be decomposed into two distinct regions: the first one, for $L=0,\dotsc,L^{\star}$ where the criterion is
growing rapidly; and the second one, for $L=L^{\star},\dotsc,L_{\text{max}}$,
where the criterion is barely increasing \citep{Lavielle:2005}. Hence, for each possible value of
$L$ in $L=1,\dotsc,L_{\text{max}}$, we compute least square linear
regressions for both parts of the graph (before and after $L$); the estimated
number of change-points is the value of $L$ that yields the best fit, that is,
the value for which the sum of the residual sums of squares computed on both
parts of the graph is minimal. For an illustration of this
methodology, see Figure~\ref{fig:rupPente}.
The case $L=0$
is treated separately and the procedure described above is used
only when the value of the test statistic $W_n$ for the presence of a
single change-point (see Section~\ref{sec:CP-onechange}) is 
significant ($p$-value smaller than 0.1\%) based on Theorem~\ref{theo:change-point}.
%\fixme{At what level?}

\begin{figure}[t] %/tsi/banda/lung/DynKW/simuCGH-fig-jrssb.py
  \centering
     \includegraphics[width=0.95\textwidth,height=2cm]{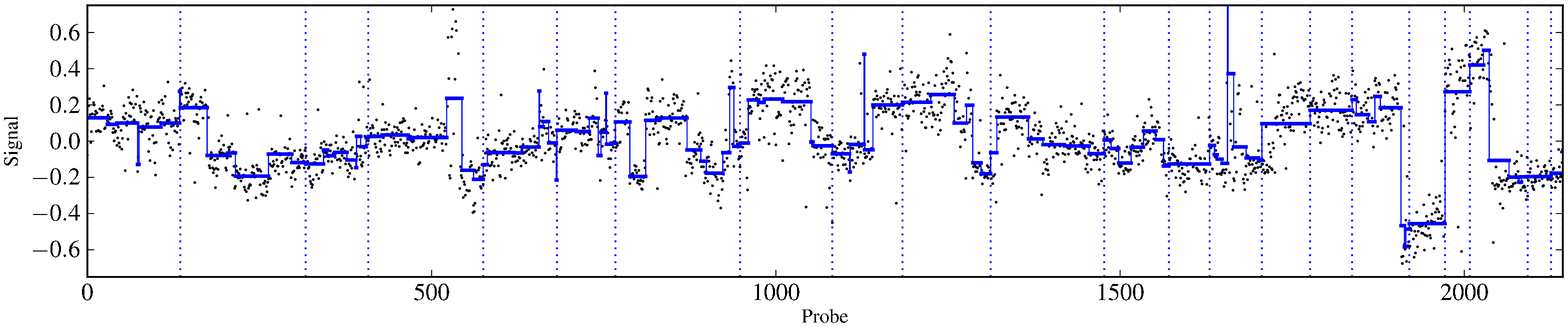}
     \includegraphics[width=0.95\textwidth,height=2cm]{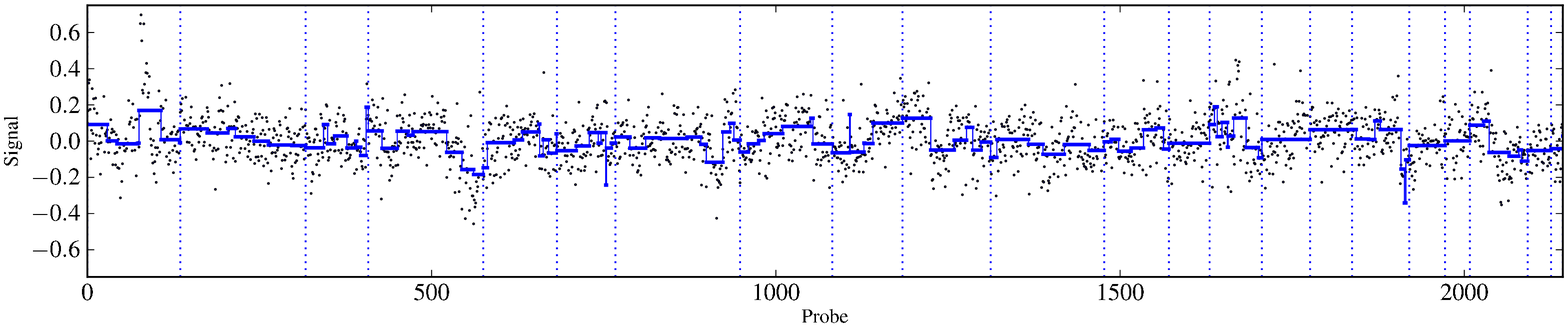}
     \includegraphics[width=0.95\textwidth,height=3cm]{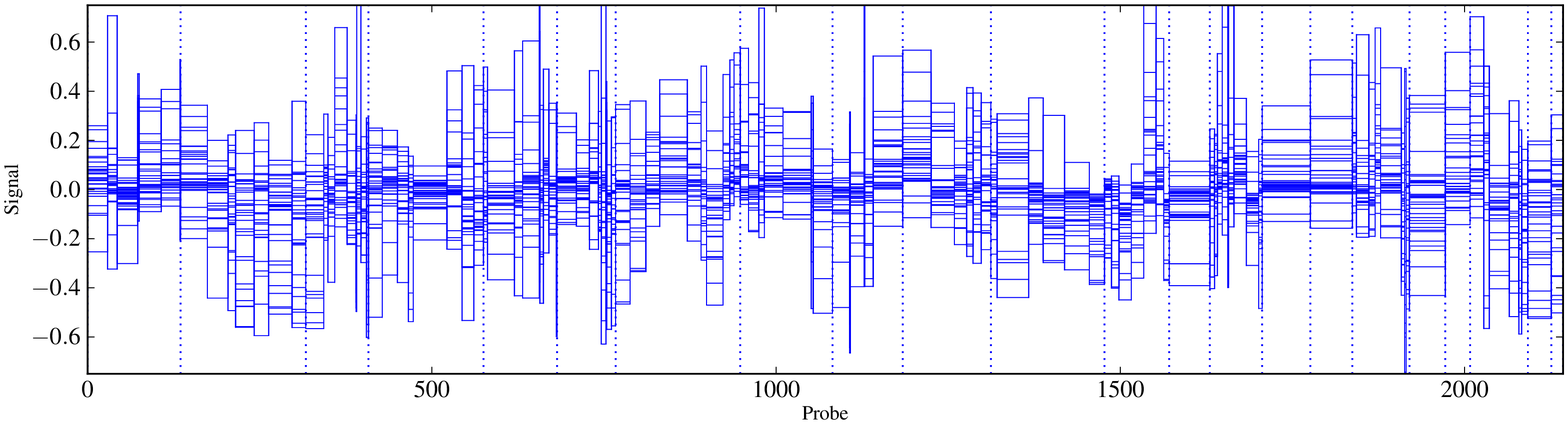}
  \caption{First and second row: copy number data for two different individuals with superimposition of the segmentation.
    Third row:  superimposition of the smoothed bladder tumour aCGH data for 32 individuals in Stage T2 cancer that result from the segmentation. 
Vertical dashed lines represent the separation between the different chromosomes.}
  \label{fig:CGH_combined}
\end{figure}

\begin{figure}[ht]
  \centering
  \includegraphics[width=.95\textwidth]{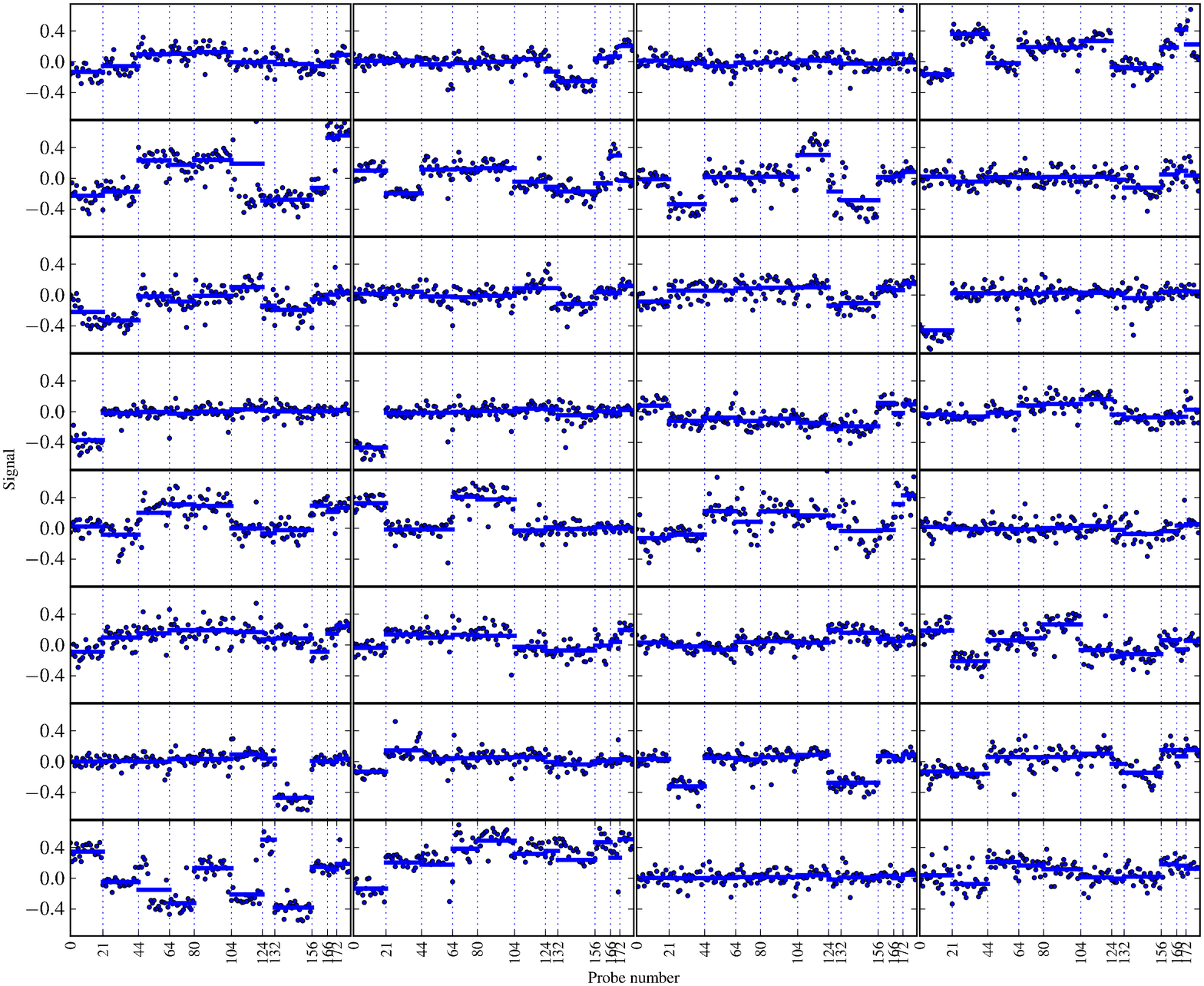}

  \caption{Data for 32 individuals in Stage T2 bladder cancer with superimposed segmentation for chromosome 7. 10 change-points were estimated and the dashed vertical lines correspond to the estimated segment boundaries.}
  \label{fig:CGH_one}
\end{figure}

Results are shown for a group of 32 profiles corresponding to Stage T2 of a tumour.
In Figure~\ref{fig:CGH_combined} the copy-number data and the segmentation of the whole set of chromosomes is displayed for two particular individuals, together with the corresponding stepwise constant approximations of the data for the 32 individuals (in the bottom of the Figure).
Figure~\ref{fig:CGH_one} details the results pertaining to the 7th chromosome.
%In Figure~\ref{fig:CGH_combined}, the smoothed version of the data of two particular individuals for the whole set of chromosomes resulting from the segmentation is displayed, while a focus is given on the 7th chromosome in Figure~\ref{fig:CGH_one}. 
In both cases, the segmentation result is represented by a signal which is constant (and equal to the mean of the data) within the detected segments.
The bottom plot in Figure~\ref{fig:CGH_combined} particularly highlights the fact that several coordinates indeed jump at the same time, suggesting that the joint segmentation model is appropriate. On the other hand, it is also obvious that one cannot assume (see, e.g., the third change-point at index 64 in Figure~\ref{fig:CGH_one}) that all coordinates undergo similar changes. Note also that in this application, the fact that the MultiRank test statistic is properly normalised with respect to the length $n$ of the data and their dimension $K$ is particularly important: $n$ corresponds to the number of probes and varies with each chromosome, $K$ represents the number of individuals and varies when considering different groups of subjects.

As a reference, the group fused Lasso algorithm by \cite{bleakley:vert;2011} outputs similar results. In particular, on the 7th chromosome, change-points are found at positions 21, 44, 65, 102, 107, 112, 124, 132, 156 and 166. On the whole set of chromosomes, 96 change-points are found while the MultiRank estimation procedure outputs 98.

\section{Conclusion}
\label{sec:conclusion}

We proposed an approach for retrospective detection of multiple changes in multivariate data. The
basic idea, used for homogeneity testing when the data groupings are known, is an extension of
well-known marginal rank based tests (Wilcoxon/Mann-Whitney and Kruskal-Wallis) based on the idea
originally proposed by \cite{wei:lachin:1984}. The use of this approach for change-point detection
(when the segments boundaries are unknown) was shown to be computationally feasible. In addition, it incorporates important aspects of the problem, in particular the fact that simultaneous detections
in different coordinates make the presence of an actual change more likely. The method was shown to
be robust against various alternatives and on a par with optimal methods in benchmark cases. The approach can also be straightforwardly modified to deal with ordinal data, missing or censored values.

To improve the method, it would be desirable to provide significance levels for the change-point
detection test when used to detect more than a single potential change-point. We believe that by
considering the normalisation used to define the statistic $W_n$ in Section~\ref{sec:CP-onechange},
it is possible to study the asymptotic behaviour of the change-point detection statistic in the
general case. This being said, the difference between the two forms of normalisation could be more
significant when applied to more than two segments. On a different level, one should obviously
consider more principled approaches for selecting the number of change-points. General-purpose
penalisation schemes could be used but we feel that novel ideas need to be developed specifically
for the change-point problem given the specific nature of over- and under-estimating  the number
of change-points. For instance, in many practical applications the significance of over-estimating
the number of change-points depends not only on the number of spurious segments but also on their
locations. Traditional approaches based on complexity penalties \citep{bai:perron:2003}, Bayesian
methods \citep{fearnhead:2006} and sparsity-based criterions
\citep{harchaoui:levy:2010,vert:bleakley:2010} are already available but there is certainly room
for new developments in these fields.

%\small{
\appendix

% \section{Appendix: Proof of (\ref{eq:ranks})}
% \label{proof:eq:ranks}

% Let $k$ be the function defined by $k(x,y)=\1(x\leq y)-\1(y\leq x)$. Note that $k$ has the following
% property: $k(x,y)=-k(y,x)$, for all $x,y$ in $\mathbb{R}$.
% Hence, we get that
% \begin{multline*}
% \sqrt{n n_1 (n-n_1)} U_{n,k}(n_1) =\sum_{i=1}^{n} \sum_{j=n_{1}+1}^{n}
% \left\{\1(X_{i,k}\le X_{j,k})-\1(X_{j,k}\le X_{i,k})\right\}\\
% =\sum_{j=n_{1}+1}^{n} \{R_j^{(k)}-\sum_{i=1}^{n} [1-\1(X_{i,k}<X_{j,k})]\}
% =\sum_{j=n_{1}+1}^{n}
% \{2R_j^{(k)}-n-\sum_{i=1}^{n}\1(X_{i,k}=X_{j,k})\}\\
% =2\sum_{j=n_{1}+1}^{n} \{R_j^{(k)}-(n+1)/2\}\;,
% \end{multline*}
% where we used for writing the last equality that there are
% no ties in the data \textit{i.e.} $X_{i,k}=X_{j,k}$ if and only if $i=j$. Using
% similar arguments, we obtain the other equality of (\ref{eq:ranks}). 

\section{Appendix: Proof of Theorem~\ref{th:WL-theorem} }
\label{sec:Appendix}

%\subsection{Proof of Theorem~\ref{th:WL-theorem}}
%\label{sec:proof1}
The proof is based on the Hoeffding decomposition of $U_{n,k}(n_1)$ for
each $k$ in $\{1,\dotsc,K\}$.
%In the sequel, for notation convenience, the index $k$ is removed.
For further details on the  Hoeffding decomposition, we refer the reader 
to Chapters  11 and 12 of \cite{vandervaart:1998}.
For each $k$ in $\{1,\dotsc,K\}$, let
 $h_{1,k}(y) = \int\! h(x,y)dF_k(x)$ and $\tilde h_{1,k}(x) = \int\! h(x,y)
dF_k(y)$, where $h$ is defined 
% for all $x,y$ in $\Bbb R$ 
by $h(x,y) = \1(x\leq y)-\1(y\leq x)$.
By the continuity of $F_k$, % we get that
$h_{1,k}(y)= 2F_k(y) - 1$ and $\tilde h_{1,k}(x)=1-2F_k(x)$.
The Hoeffding decomposition of $U_{n,k}(n_1)$ can thus be written as
$
U_{n,k}(n_1)=\hat U_{n,k}(n_1)+R_{n,k}(n_1)
$, 
where
\begin{align}\label{eq:U_nHoeffding_dec1}
  & \hat U_{n,k}(n_1) = \frac{n_{1}}{\sqrt{n n_1 (n-n_1)}} \sum_{j=n_{1}+1}^{n}h_{1,k}(X_{j,k}) +
  \frac{n-n_{1}}{\sqrt{n n_1 (n-n_1)}} \sum_{i=1}^{n_1}\tilde h_{1,k}(X_{i,k})\; , \\
\label{eq:U_nHoeffding_dec2}
 & R_{n,k}(n_1)=\inv{\sqrt{n n_1 (n-n_1)}}
  \sum_{i=1}^{n_{1}}\sum_{j=n_{1}+1}^{n}\![h(X_{i,k},X_{j,k})-\tilde h_{1,k}(X_{i,k})-h_{1,k}(X_{j,k})]\; .
\end{align}

We first prove that $U_{n,k}(n_1)=\hat U_{n,k}(n_1)+o_p(1)$ by showing that $\var[R_{n,k}(n_1)]\to 0$, as $n$ tends to infinity.
Using that $\PE[U_{n,k}(n_1)]=\PE[\hat U_{n,k}(n_1)]=0$, we obtain that
$\var[R_{n,k}(n_1)]=\var[U_{n,k}(n_1)-\hat U_{n,k}(n_1)]
%=\var[U_{n,k}(n_1)]+\var[\hat U_{n,k}(n_1)]-2\cov[U_{n,k}(n_1)\hat U_{n,k}(n_1)]
=\PE[U_{n,k}^2(n_1)]+\PE[\hat U_{n,k}^2(n_1)]-2\PE[U_{n,k}(n_1)\hat U_{n,k}(n_1)].$
By independence of the $(X_{i,k})_{1\leq i\leq n}$, we obtain that
\begin{multline}\label{eq:varUn_hat_bis}
  \PE[\hat U_{n,k}^2(n_1)] = \frac{n_1^2}{n n_1 (n-n_1)}\sum_{j=n_1+1}^n \PE[h_{1,k}(X_{j,k})^2] +
  \frac{(n-n_1)^2}{n n_1 (n-n_1)} \sum_{i=1}^{n_1}\PE[\tilde
  h_{1,k}(X_{i,k})^2].
\end{multline}
Using that
\begin{equation}
  \label{eq:variance_h}
  \PE[h_{1,k}(X_{i,k})^2] = 4\PE[(F_k(X_{1,k})-1/2)^2] =
  4\var(\mathcal{U}) = 1/3\;,
\end{equation}
where $\mathcal{U}$ has a uniform distribution on $[0,1]$, we get, on
the one hand, that
\begin{equation}\label{eq:varUn_hat}
  \PE[\hat U_{n,k}^2(n_1)] = \frac{n_1^2(n-n_1)}{3 n n_1 (n-n_1)} + \frac{ (n-n_1)^2n_1}{3 n n_1 (n-n_1)}=1/3\;.
\end{equation}
On the other hand
\begin{multline} \label{eq:varUn}
\PE[U_{n,k}^2(n_1)] % = \inv{n n_1 (n-n_1)}\sum_{1\leq i,i'\leq n_1} \sum_{n_1+1\leq j,j'\leq n}
%   \PE[h(X_{i,k},X_{j,k})h(X_{i',k},X_{j',k})]\\
%  \label{eq:varUn_1}
  = \inv{n n_1 (n-n_1)} \sum_{i=1}^{n_{1}}\sum_{j=n_{1}+1}^{n} \PE[h(X_{i,k},X_{j,k})^2] \\
%  \label{eq:varUn_2}
  + \inv{n n_1 (n-n_1)} \sum_{1\leq i\neq i'\leq n_1} \sum_{j=n_1+1}^n
  \PE[h(X_{i,k},X_{j,k})h(X_{i',k},X_{j,k})]\\
%  \label{eq:varUn_3}
  + \inv{n n_1 (n-n_1)} \sum_{i=1}^{n_1} \sum_{n_1+1\leq j\neq j'\leq n}\PE[h(X_{i,k},X_{j,k})h(X_{i,k},X_{j',k})] \; .
\end{multline}
We separately study the three terms of the r.h.s of (\ref{eq:varUn}).
% (\ref{eq:varUn_2}) and (\ref{eq:varUn_3}) of equation
% (\ref{eq:varUn}).
Using that $(X_{i,k})_{1\leq i\leq n}$ are i.i.d. , we get 
\begin{equation}\label{eq:varUn_1Calc}
\inv{n n_1 (n-n_1)} \sum_{i=1}^{n_{1}}\sum_{j=n_{1}+1}^{n} \PE[h(X_{i,k},X_{j,k})^2]
= \frac{n_1(n-n_1)}{n n_1 (n-n_1)} \PE[h(X_{1,k},X_{n_1+1,k})^2]\to 0,\; \text{as  $n\to\infty$} \; .
\end{equation}
%which converges to 0 as $n$ goes to infinity. 
Then, by continuity of $F_k$, we have
\begin{multline}\label{eq:varUn_2Calc}
\inv{n n_1 (n-n_1)} \sum_{1\leq i\neq i'\leq n_1} \sum_{j=n_1+1}^n
  \PE[h(X_{i,k},X_{j,k})h(X_{i',k},X_{j,k})]\\
%  = \frac{n_1^2(n-n_1)}{n n_1 (n-n_1)} \iiint  (\1(x\leq y)-\1(y\leq x))(\1(z\leq y)-\1(y\leq z)) dF_k(x)dF_k(y)dF_k(z)\\
  =  \frac{(n_1^2-n_1)(n-n_1)}{n n_1 (n-n_1)} \int (2F_k(y)-1)(2F_k(y)-1)dF_k(y)
  =  \frac{n_1(n_1-1)(n-n_1)}{3n n_1 (n-n_1)}\; .
\end{multline}
Using similar arguments, the last term of the r.h.s of (\ref{eq:varUn}) is equal
to $n_1(n-n_1)(n-n_1-1)/(3n n_1 (n-n_1))$. With \eqref{eq:varUn_1Calc} and \eqref{eq:varUn_2Calc},
we obtain
\begin{equation}\label{eq:varUn_3Calc}
\PE[U_{n,k}^2(n_1)]\to 1/3,\; \text{as $n\to\infty$}.
\end{equation}
% \begin{equation}
%   \label{eq:varUn_3Calc}
%   (\ref{eq:varUn_3}) =  \frac{n_1(n-n_1)^2}{3n^3}
% \end{equation}
% Using \eqref{eq:varUn_1Calc}, \eqref{eq:varUn_2Calc} and
% \eqref{eq:varUn_3Calc}, we can thus conclude that
% \begin{equation}
%   \label{eq:varUnCls}
%   \var(U_n) \stackrel{n\rightarrow \infty}{\longrightarrow}  \frac{n_1(1-n_1)}{3}
% \end{equation}
Since $\PE[U_{n,k}(n_1)\hat U_{n,k}(n_1)]\to 1/3$,
as $n\to\infty$, \eqref{eq:varUn_hat} and \eqref{eq:varUn_3Calc} lead to
$\var[R_{n,k}(n_1)]\to 0$ and thus $U_{n,k}(n_1)=\hat U_{n,k}(n_1)+o_p(1)$,
as $n$ tends to infinity.
The multivariate central limit theorem then yields
$
(U_{n,1}(n_1),\dotsc,U_{n,K}(n_1))'\to\mathcal{N}(0,\Sigma)\; ,
$
where the $(k,k')$th entry of $\Sigma$ is given by 
$\Sigma_{kk'} =\lim_{n\to\infty}\PE[\hat U_{n,k}(n_1)\hat U_{n,k'}(n_1)]$.
Using that the $(X_{i,k})_{1\leq i\leq n}$ are i.i.d., we obtain that
\begin{multline*}
\PE[\hat U_{n,k}(n_1)\hat U_{n,k'}(n_1)]
=\frac{4n_1^2}{n n_1 (n-n_1)}\sum_{j=n_1+1}^n \PE[\{F_k(X_{j,k})-1/2\}\{F_{k'}(X_{j,k'})-1/2\}]\\
+\frac{4(n-n_1)^2}{n n_1 (n-n_1)}\sum_{i=1}^{n_1} \PE[\{F_k(X_{i,k})-1/2\}\{F_{k'}(X_{i,k'})-1/2\}] \\
%=\left(\frac{4n_1}{n}+\frac{4(n-n_1)}{n}\right)\PE[\{F_k(X_{1,k})-1/2\}\{F_{k'}(X_{1,k'})-1/2\}]
%=\frac{4n_1(n-n_1)}{n n_1 (n-n_1)}\PE[\{F_k(X_{1,k})-1/2\}\{F_{k'}(X_{1,k'})-1/2\}]\\
= 4\cov\left(F_{k}(X_{1,k}),F_{k'}(X_{1,k'})\right)\; .
\end{multline*}
Thus,
$\Sigma^{-1/2}(U_{n,1}(n_1),\dotsc,U_{n,K}(n_1))'\inlaw\mathcal{N}(0,\operatorname{Id}_K)$. Since
$\hat{\Sigma}_n\inproba\Sigma$, we deduce from Slutsky's Theorem
that
$\hat{\Sigma}_n^{-1/2}(U_{n,1}(n_1),\dotsc,U_{n,K}(n_1))'\inlaw\mathcal{N}(0,\operatorname{Id}_K)$,
which concludes the proof.

\section{Appendix: Proof of Theorem~\ref{th:MKW}}
\label{sec:proofMKW}

Using that $R_j^{(k)}=\sum_{i=1}^n \1(X_{i,k}\leq X_{j,k})$, we obtain that
\begin{multline}\label{eq:R_bullet}
\bar{R}_{\ell}^{(k)}-\frac{n+1}{2}=\frac{1}{n_{\ell+1}-n_{\ell}}\left(\sum_{j=n_{\ell}+1}^{n_{\ell+1}}
\sum_{i=1}^n \1(X_{i,k}\leq X_{j,k})\right)-\frac{n+1}{2}\\
=\frac{1}{n_{\ell+1}-n_{\ell}}\sum_{j=n_{\ell}+1}^{n_{\ell+1}}
\sum_{\stackrel{i=1}{i\neq j}}^n \left[\1(X_{i,k}\leq X_{j,k})-1/2\right]\;.
\end{multline}
Let $h(x,y)=\1(x\leq y)$, $h_{1,k}(y)=\int\1(x\leq
  y)\rmd F_{k}(x)$ and $h_{2,k}(x)=\int\1(x\leq y)\rmd F_{k}(y).$
By continuity of $F_{k}$:
$h_{1,k}(y)=F_{k}(y)$ and  $h_{2,k}(x)=1-F_{k}(x).$ Using the notation:
$$
\mathcal{R}_{\ell}^{(k)}=(n_{\ell+1}-n_{\ell})^{1/2}/n(\bar{R}_{\ell}^{(k)}-(n+1)/2)\;,
$$
the Hoeffding decomposition yields
\begin{multline}\label{eq:dec:hoeffding}
\mathcal{R}_{\ell}^{(k)}
=\frac{(n_{\ell+1}-n_{\ell})^{1/2}}{n}\left[\frac{n-1}{n_{\ell+1}-n_{\ell}}\sum_{j=n_{\ell}+1}^{n_{\ell+1}}(h_{1,k}(X_{j,k})-1/2)
%\right]
%\\
%+\frac{(n_{\ell+1}-n_{\ell})^{1/2}}{n}\left[
+\frac{n_{\ell+1}-n_{\ell}-1}{n_{\ell+1}-n_{\ell}}\sum_{i=1}^{n}(h_{2,k}(X_{i,k})-1/2)\right]\\
+\frac{(n_{\ell+1}-n_{\ell})^{1/2}}{n}\left[\frac{1}{n_{\ell+1}-n_{\ell}}\sum_{j=n_{\ell}+1}^{n_{\ell+1}}
\sum_{\stackrel{i=1}{i\neq j}}^n \left\{h(X_{i,k},X_{j,k})-h_{1,k}(X_{j,k})-h_{2,k}(X_{i,k})+1/2\right\}\right]\\
\stackrel{\textrm{def}}{=}\mathcal{R}_{\ell,1}^{(k)}+\mathcal{R}_{\ell,2}^{(k)}+\mathcal{R}_{\ell,3}^{(k)}\;.
\end{multline}
Note that $\mathcal{R}_{\ell,3}^{(k)}=o_p(1)$, as $n$ tends to infinity, since it can be proved that
$\var(\mathcal{R}_{\ell,3}^{(k)})=\var[\mathcal{R}_{\ell}^{(k)}-(\mathcal{R}_{\ell,1}^{(k)}+\mathcal{R}_{\ell,2}^{(k)})]\to
0$, as $n$ tends to infinity.
Thus,  (\ref{eq:dec:hoeffding}) can be rewritten as
\begin{multline*}
\mathcal{R}_{\ell}^{(k)}=\frac{n-1}{n(n_{\ell+1}-n_{\ell})^{1/2}}\sum_{j=n_{\ell}+1}^{n_{\ell+1}}(F_{k}(X_{j,k})-1/2)
+\frac{n_{\ell+1}-n_{\ell}-1}{n(n_{\ell+1}-n_{\ell})^{1/2}}\sum_{i=1}^{n}(1/2-F_{k}(X_{i,k}))
+o_p(1)\;. 
\end{multline*}
Since $\sum_{i=1}^{n}(1/2-F_{k}(X_{i,k}))=\sum_{p=0}^{L-1}\sum_{j=n_p+1}^{n_{p+1}}(1/2-F_{k}(X_{j,k}))$,
\begin{multline*}
\mathcal{R}_{\ell}^{(k)}=\frac{n-(n_{\ell+1}-n_{\ell})}{n(n_{\ell+1}-n_{\ell})^{1/2}}\sum_{j=n_{\ell}+1}^{n_{\ell+1}}(F_{k}(X_{j,k})-1/2)
-\frac{n_{\ell+1}-n_{\ell}-1}{n(n_{\ell+1}-n_{\ell})^{1/2}}\sum_{\stackrel{p=0}{p\neq\ell}}^{L-1}\sum_{j=n_p+1}^{n_{p+1}}(F_{k}(X_{j,k}-1/2))+o_p(1)\\
\stackrel{\textrm{def}}{=}U_k(n_\ell,n_{\ell+1})+o_p(1)\;. 
\end{multline*}
Observe that, for a fixed $\ell$ in $\{0,\dots,L-1\}$ and
$k,k'$ in $\{1,\dots,K\}$, we get, as $n$ tends to infinity,
\begin{multline}\label{eq:lim1}
4\cov(U_k(n_\ell,n_{\ell+1}),U_{k'}(n_\ell,n_{\ell+1}))=\\
\Sigma_{kk'}\left[\left(1-\frac{(n_{\ell+1}-n_\ell)}{n}\right)^2
+\sum_{\stackrel{p=0}{p\neq
    \ell}}^{L-1}\frac{(n_{\ell+1}-n_{\ell}-1)^2(n_{p+1}-n_p)}{n^2(n_{\ell+1}-n_{\ell})}\right]
\to (1-t_{\ell+1})\Sigma_{kk'}\;,
\end{multline}
where we have used that $\sum_{\stackrel{p=0}{p\neq\ell}}^{L-1}(n_{p+1}-n_p)=n-(n_{\ell+1}-n_{\ell}).$
In the same way, for fixed $k,k'$ in $\{1,\dots,K\}$ and $\ell\neq
\ell'$ in $\{0,\dots,L-1\}$, we obtain,  as $n$ tends to
infinity,
\begin{equation}\label{eq:lim2}
4\cov(U_k(n_\ell,n_{\ell+1}),U_{k'}(n_{\ell'},n_{\ell'+1}))\to -\sqrt{t_{\ell+1}t_{\ell'+1}}\Sigma_{kk'}\;.
\end{equation}
Let 
$$
\bar{R}_n=2\left(\frac{(n_{1}-n_0)^{1/2}}{n}\bar{R}_0',\dots,\frac{(n_{L}-n_{L-1})^{1/2}}{n}\bar{R}_{L-1}'\right)'\;.
$$
We deduce from (\ref{eq:lim1}), (\ref{eq:lim2}) and the multivariate
central limit  theorem that
$$
\bar{R}_n
\stackrel{d}{\longrightarrow}\mathcal{N}(0,\Theta\otimes\Sigma)\;, n\to\infty\;,
$$
where $\Sigma$ is the $K\times K$ matrix defined in
(\ref{eq:CovMatrix}), $\otimes$ denotes the Kronecker product,
$\Theta=\operatorname{Id}_{L}-\boldsymbol{\sqrt{t}\sqrt{t}}'$
with $\boldsymbol{\sqrt{t}}=(\sqrt{t_1},\dots,\sqrt{t_{L}})'$.
Thus,
$$
%\frac{(n_{k+1}-n_k)^{1/2}}{n}\Sigma^{-1/2}(\bar{R}_k^{(\ell)}-n/2)
\bar{R}_n^{\Sigma}
\stackrel{d}{\longrightarrow}\mathcal{N}(0,\Theta\otimes\operatorname{Id}_{K})\;, n\to\infty\;,
$$
where 
$$
\bar{R}_n^{\Sigma}=2\left(\frac{(n_{1}-n_0)^{1/2}}{n}\Sigma^{-1/2}\bar{R}_0',\dots,\frac{(n_{L}-n_{L-1})^{1/2}}{n}\Sigma^{-1/2}\bar{R}_{L-1}'\right)'\;.
$$
Since $\hat{\Sigma}_n\stackrel{p}{\longrightarrow}\Sigma$, as $n$
tends to infinity, the same convergence
holds when $\Sigma$ is replaced by $\hat{\Sigma}_n$.
Since $\sum_{\ell=0}^{L-1}(n_{\ell+1}-n_\ell)/n=1$,
$\sum_{\ell=1}^{L}t_\ell=1$ and the matrix $t$ has
eigenvalue 0 of multiplicity 1 (with eigenspace spanned by
$\boldsymbol{\sqrt{t}}$), and eigenvalue 1 of multiplicity
$L-1$. Hence, the eigenvalues of $\Theta\otimes\operatorname{Id}_{K}$
are 0, with multiplicity $K$, and 1, with multiplicity
$(L-1)K$, which concludes the proof using Cochran's theorem.

\section{Appendix: Proof of Theorem \ref{th:power}}
\label{appendix:power}
In Appendix \ref{sec:Appendix}, we proved that under the null hypothesis
where the $\mathbf{X}_j$'s are i.i.d. random vectors such that the
c.d.f $F_k$ of $X_{1,k}$ is continuous:
\begin{multline}\label{eq:Unk}
U_{n,k}(n_1) = \frac{2n_{1}}{\sqrt{n n_1 (n-n_1)}} \sum_{j=n_{1}+1}^{n}\{F_{k}(X_{j,k})-1/2\} \\
-\frac{2(n-n_{1})}{\sqrt{n n_1 (n-n_1)}} \sum_{i=1}^{n_1}\{F_{k}(X_{i,k})-1/2\}+o_P(1)\;, 
\textrm{ as } n\to\infty\;.
\end{multline}
Since, by assumption, the model $f(\mathbf{x}-\boldsymbol{\theta})$ is differentiable in quadratic mean at
$\theta$, the log-likelihood ratio $L_n$ defined by
$$
L_n=\log\left[\frac{\prod_{i=1}^{n_1} f(\mathbf{X}_i)\prod_{j=n_1+1}^n f(\mathbf{X}_j-\boldsymbol{\delta}/\sqrt{n})}
{\prod_{i=1}^{n} f(\mathbf{X}_i)}\right]
$$
satisfies the following asymptotic expansion as $n$ tends to infinity
\begin{equation}\label{eq:Ln}
L_n=-\frac{\boldsymbol{\delta}}{\sqrt{n}}\sum_{j=n_1+1}^n \nabla\log f(\mathbf{X}_j)
-\frac {(1-t_1)}{2} \boldsymbol{\delta}' I_f \boldsymbol{\delta}+o_P(1)\;.
\end{equation}
Combining multivariate central limit theorem with (\ref{eq:Unk}) and (\ref{eq:Ln}) one can thus show that $(\mathbf{U}_n,L_n)$ converges in distribution to a Gaussian random vector. From the expression of the asymptotic covariance matrix and
using Le Cam's third lemma, one obtains
that under $(H_{1,n})$
$$
\mathbf{U}_n(n_1)\stackrel{d}{\longrightarrow}
\mathcal{N}_K(2\sqrt{t_1(1-t_1)} A\boldsymbol{\delta},\Sigma)\;.
$$
Using Glivenko-Cantelli Theorem, the weak law of large numbers and the
fact that $f_k$ is bounded for all $k$ --which implies that $F_k$ is a
Lipschitz function for all $k$-- it can be proved that, under $(H_{1,n})$, $\hat{\Sigma}_n$
converges in probability to $\Sigma$. The conclusion follows
using Slutsky's Lemma.
%The convergence in (\ref{eq:conv_hotelling}) is obtained by using the
%multivariate central limit theorem.
\section{Appendix: Proof of Corollary  \ref{coro:indep}}
\label{appendix:coro:indep}
Assuming independence, $\Sigma = 4 \PE_0[\{F_1(X_{1,1})-1/2\}^2] \operatorname{Id}_K = 1/3 \operatorname{Id}_K$, 
where $\operatorname{Id}_K$ denotes the $K\times K$
 identity matrix, and $A$ is the diagonal matrix with elements
$A_{k,k}=\int_{\mathbb{R}} F_k(x) f_k'(x)\rmd x=-\int_{\mathbb{R}} f_k^2(x)\rmd x
=-\sigma_k^{-1}\int_{\mathbb{R}} (\sigma_k f_k(\sigma_k x))^2\rmd x$. 
The lower bound on $d_S$
is obtained by the classical result that $12\lambda_k^2\geq 108/125$ \cite[p. 198]{vandervaart:1998}. Finally, in the Gaussian case, $\lambda_k=1/(2\sqrt{\pi})$.

\section{Appendix: Proof of Corollary \ref{coro:gaussian}}
\label{appendix:coro:gaussian}
Let us first prove that $A$ is a diagonal matrix such that
$A_{k,k}=\sigma_k^{-1}/(2\sqrt{\pi})$. Let $D=(d_{k,\ell})_{1\leq k,\ell\leq K}=C^{-1}$, then
$$
A_{k,\ell}=\PE_0[(F_k(X_{1,k})-1/2)(\sum_{j=1}^K d_{\ell,j} X_{1,j})]
=\sum_{j=1}^K d_{\ell,j} \sigma_j \PE[\Phi(\tilde{X}_{1,k})\tilde{X}_{1,j}]\;,
$$
where $\Phi$ is the c.d.f. of a standard Gaussian random variable
and $\tilde{X}_{1,j}=X_{1,j}/\sigma_j$ is a standard Gaussian random variable.
Since $\Phi(x)-1/2=\sum_{p\geq 1}(\alpha_p/p!)H_p(x)$, where $\alpha_p=\PE[\Phi(Z)H_p(Z)]$,
$Z$ is a standard Gaussian random variable and $H_p$ is the $p$th
Hermite polynomial with leading coefficient equal to one ($H_1(x)=x$,
$H_2(x)=x^2-1$, ...). By applying Mehler's formula, one obtains 
$\PE[\Phi(\tilde{X}_{1,k})\tilde{X}_{1,j}]=\alpha_1
C_{k,j}/(\sigma_k\sigma_j)$, where
$\alpha_1=\PE[Z\Phi(Z)]=\int_{\mathbb{R}}\varphi^2(x)\rmd x$, using an
integration by parts.
Thus,
$A_{k,\ell}=\sigma_k^{-1}/(2\sqrt{\pi})\sum_{j=1}^K d_{\ell,j} C_{j,k}=\sigma_k^{-1}/(2\sqrt{\pi})\1_{\{k=\ell\}}.$
The expression given for $\Sigma$ can be obtained similarly and is well-known in the literature as it provides the link between the Spearman correlation and the usual correlation coefficients~\citep{kruskal:1958}. Regarding the lower bound, first note that
\begin{multline*}
\textrm{ARE}=(4\pi)^{-1}\boldsymbol{\delta}'\textrm{diag}(\sigma_1^{-1},\dots,\sigma_K^{-1}) (4\Sigma^{-1})
\textrm{diag}(\sigma_1^{-1},\dots,\sigma_K^{-1})\boldsymbol{\delta}/\boldsymbol{\delta}' C^{-1} \boldsymbol{\delta}\\
\geq (4\pi)^{-1}\lambda_{\min}\left(\textrm{diag}(\sigma_1^{-1},\dots,\sigma_K^{-1}) (4\Sigma^{-1})
\textrm{diag}(\sigma_1^{-1},\dots,\sigma_K^{-1}) C\right)\\
=(4\pi)^{-1}\sigma_{\max}^{-2}\lambda_{\min}(C)/\lambda_{\max}(\Sigma/4)\; .
\end{multline*}
From the expression of $\Sigma$, it is easily checked that $\Sigma_{k,\ell} \leq |C_{k,\ell}/3\sigma_k^{-1}\sigma_\ell^{-1}| \leq |C_{k,\ell}|/3\sigma_{\min}^{-2}$,
which gives the second result.

\section{Appendix: Proof of Theorem~\ref{theo:change-point}}
\label{sec:Appendix2}

We start by proving (\ref{eq:WLrupture}) when $\hat{\Sigma}_n$ is
replaced by $\Sigma$ in (\ref{eq:S_n_tilde}).
For this, we shall verify the assumptions of \cite[Theorem
15.6]{billingsley:1968}: the convergence of the finite-dimensional
distributions:
\begin{multline}
  \label{eq:fini-dim}
  \left( \mathbf{V}_{n}(\lfloor n t_1\rfloor)'\Sigma^{-1}
    \mathbf{V}_{n}(\lfloor n t_1\rfloor),\dotsc,
    \mathbf{V}_{n}(\lfloor n t_p\rfloor)'\Sigma^{-1}
    \mathbf{V}_{n}(\lfloor n t_p\rfloor) \right) \\
  \inlaw \left(
    \sum_{k=1}^{K} B_{k}^{2}(t_1),\dotsc,\sum_{k=1}^{K} B_{k}^{2}(t_p)
  \right), \textrm{ for } 0<t_1<\dotsc<t_p<1\;,\; n\to\infty\;,
\end{multline}
and the tightness criterion for the process:
$$
\left\{
  \mathbf{V}_{n}(\lfloor n t\rfloor)'\Sigma_n^{-1} \mathbf{V}_{n}(\lfloor n t\rfloor)
;\; 0<t<1
\right\}\;,
$$
where $\lfloor x\rfloor$ denotes the integer part of $x$.
Let $n_1=\lfloor n t_1\rfloor$, with $t_1$ in $(0,1)$.
In the same way as in Appendix~\ref{sec:Appendix}, 
as $V_{n,k}(\cdot)$ only differs from $U_{n,k}(\cdot)$ by a
normalising factor, we can prove that
$V_{n,k}(n_1)=\hat V_{n,k}(n_1)+o_p(1)$, with $0<n_1<n$ and $$\hat
V_{n,k}(n_1) = \frac{n_1}{n^{3/2}} \sum_{j=n_1+1}^n
h_{1,k}(X_{j,k}) - \frac{n-n_1}{n^{3/2}}\sum_{i=1}^{n_1} h_{1,k}(X_{i,k})\;,$$
where $h_{1,k}(x)=2F_k(x)-1$ and that 
\begin{equation}
  \label{eq:covVHat}
  \PE[\hat V_{n,k}(n_1) \hat V_{n,k'}(n_1)] \to
  4t_1(1-t_1)\cov\left(F_{k}(X_{1,k}),F_{k'}(X_{1,k'})\right),\;
  \text{as $n\to\infty$} \; . 
\end{equation}
Let $n_2=\lfloor n t_2\rfloor$. Since $1<n_1<n_2<n$, $n_1/n \to t_1 \in (0,1)$, and $n_2/n \to t_2 \in (0,1)$  we get
\begin{multline}
  \label{eq:EspUkn1k'n2}
  \PE [\hat V_{n,k}(n_1) \hat V_{n,k'}(n_2)] = 
  \PE \left[ 
    \left\{
        \frac{n_1}{n^{3/2}}\sum_{j={n_1+1}}^nh_{1,k}(X_{j,k}) -
        \frac{n-n_1}{n^{3/2}}\sum_{i=1}^{n_1}h_{1,k}(X_{i,k})
      \right\} 
      \right . \\
      \times
      \left .
    \left\{
        \frac{n_2}{n^{3/2}}\sum_{j={n_2+1}}^nh_{1,k'}(X_{j,k'}) -
        \frac{n-n_2}{n^{3/2}}\sum_{i=1}^{n_2}h_{1,k'}(X_{i,k'})
      \right\}
  \right] \;.
\end{multline}
By decomposing the interval $[n_1+1,n]$ (resp. $[1,n_2]$) into
$[n_1+1,n_2]$ and $[n_2+1,n]$ (resp. $[1,n_1]$ and $[n_1+1,n_2]$)
and developing the expression, we obtain 
\begin{multline}
  \label{eq:EspUkn1k'n2_bis}
  \PE [\hat V_{n,k}(n_1) \hat V_{n,k'}(n_2)] = \\
  \PE
  \left [
    \frac{(n-n_1)(n-n_2)}{n^3} \sum_{i=1}^{n_1}h_{1,k}(X_{i,k})h_{1,k'}(X_{i,k'})
    -\frac{n_1(n-n_2)}{n^3}
    \sum_{j=n_1+1}^{n_2}h_{1,k}(X_{i,k})h_{1,k'}(X_{i,k'}) \right. \\
\left .    +\frac{n_1 n_2}{n^3} \sum_{j=n_2+1}^{n}h_{1,k}(X_{i,k})h_{1,k'}(X_{i,k'})
  \right ] \\
  = \frac{n_1(n-n_2)}{n^2} \Sigma_{kk'} 
  \longrightarrow t_1(1-t_2) \Sigma_{kk'}, \; \text{as $n\to\infty$} \; . 
\end{multline}
With \eqref{eq:covVHat} and \eqref{eq:EspUkn1k'n2_bis}, we obtain
\begin{equation}
  \label{eq:fini-dim_blocmat}
  \begin{pmatrix}
    \hat{\mathbf{V}}_{n}(n_1)\\
    \hat{\mathbf{V}}_{n}(n_2)
  \end{pmatrix}
\inlaw
\mathcal{N}\left(0;
  \left(
  \begin{array}{c|c}
    t_1(1-t_1) \Sigma &  t_1(1-t_2) \Sigma\\
    \hline
    t_1(1-t_2) \Sigma &  t_2(1-t_2) \Sigma 
  \end{array}
  \right)
\right),
\end{equation}
which is equivalent to
\begin{equation}
  \label{eq:equiLFD}
    \begin{pmatrix}
    \hat{\mathbf{V}}_{n}(n_1)\\
    \hat{\mathbf{V}}_{n}(n_2)
  \end{pmatrix}
\inlaw
\begin{pmatrix}
  \Sigma^{\inv{2}} & 0\\
  0 & \Sigma^{\inv{2}}
\end{pmatrix}
\begin{pmatrix}
\mathbf{B}(t_1)\\
\mathbf{B}(t_2)  
\end{pmatrix},
\end{equation}
where $\mathbf{B}(t)=(B_1(t),\dotsc,B_K(t))$, $0\leq t\leq 1$.
For the sake of clarity and without loss of generality,
\eqref{eq:fini-dim} is thus proved  in the case $p=2$
by applying the continuous function 
\begin{equation}
  \label{eq:apply}
  \begin{pmatrix}
    x_1\\x_2
  \end{pmatrix}
  \longmapsto
  \begin{pmatrix}
    x_1'x_1\\
    x_2'x_2
  \end{pmatrix}, \;\mathrm{where}\ x_1,x_2\in \mathbb{R}^K.
\end{equation}
% \begin{equation}
%   \begin{pmatrix}
%     \Sigma^{-1/2}  \mathbf{V}_{n}(\lfloor n t_1\rfloor)\\
%     \hline
%     \Sigma^{-1/2}  \mathbf{V}_{n}(\lfloor n t_2\rfloor)
%   \end{pmatrix}
%   \inlaw
%   \begin{pmatrix}
%     B_1(t_1)\\
%     \vdots\\
%     B_K(t_1)\\
%     \hline
%     B_1(t_2)\\
%     \vdots\\
%     B_K(t_2)
%   \end{pmatrix}
% \end{equation}
% that is
In the following, we check the tightness condition, that is, for 
$0< t_1 < t < t_2 <1$, we show that
\begin{multline}
  \label{eq:tightnessCondition}
  \PE
  \left[
  |\hat{\mathbf{V}}_n(\lfloor n t \rfloor)\Sigma^{-1}
  \hat{\mathbf{V}}_n(\lfloor n t \rfloor)
  - \hat{\mathbf{V}}_n(\lfloor n t_1 \rfloor)\Sigma^{-1}
  \hat{\mathbf{V}}_n(\lfloor n t_1 \rfloor)|^2 \right . \\
\left . \times
  |\hat{\mathbf{V}}_n(\lfloor n t_2 \rfloor)\Sigma^{-1}
  \hat{\mathbf{V}}_n(\lfloor n t_2 \rfloor)
  - \hat{\mathbf{V}}_n(\lfloor n t \rfloor)\Sigma^{-1}
  \hat{\mathbf{V}}_n(\lfloor n t \rfloor)|^2 
\right]
\leq C |t_2-t_1|^2,
\end{multline}
where $C$ is a positive constant.
Let $\mathbf{x}_n(t) =
(x_{n,1}(t),\dots,x_{n,K}(t))' = A \hat{\mathbf{V}}_n(\lfloor n t
\rfloor)$, where $A = \Sigma^{-\inv{2}}$, whose $(p,q)$th entry is
denoted by $a_{p,q}$.
The l.h.s.\ of \eqref{eq:tightnessCondition} can thus be
rewritten as
\begin{multline}
  \label{eq:tightnessCondition_bis}
  \PE
  \left[
    |\mathbf{x}'_n(t)\mathbf{x}_n(t)-\mathbf{x}'_n(t_1)\mathbf{x}_n(t_1)|^2
    |\mathbf{x}'_n(t_2)\mathbf{x}_n(t_2)-\mathbf{x}'_n(t)\mathbf{x}_n(t)|^2
  \right]\\
  = \PE
  \left[
    \left|\sum_{k=1}^K \{x_{n,k}^2(t)-x_{n,k}^2(t_1)\}\right|^2
    \left|\sum_{k'=1}^K \{x_{n,k'}^2(t_2)-x_{n,k'}^2(t)\}\right|^2
  \right].
\end{multline}
Note that
\begin{multline}
  \label{eq:tight_1}
  \sum_{k=1}^K \{x_{n,k}^2(t)-x_{n,k}^2(t_1)\} =
  \sum_{k=1}^K (x_{n,k}(t)-x_{n,k}(t_1))(x_{n,k}(t)+x_{n,k}(t_1))\\
  = \sum_{k=1}^K 
  \left(
    \sum_{p=1}^K a_{k,p}[\hat{V}_{n,p}(\floor{n t})-\hat{V}_{n,p}(\floor{n t_1})]
  \right)
  \left(
    \sum_{p'=1}^K a_{k,p'}[\hat{V}_{n,p'}(\floor{n t})+\hat{V}_{n,p'}(\floor{n t_1})]
  \right)\\
  = \sum_{p,p'=1}^K b_{p,p'}
\left[\hat{V}_{n,p}(\floor{n t})-\hat{V}_{n,p}(\floor{n t_1})\right]
\left[\hat{V}_{n,p'}(\floor{n t})+\hat{V}_{n,p'}(\floor{n t_1})\right],
  \end{multline}
  where $b_{p,p'}=\sum_{k=1}^K a_{k,p}a_{k,p'}$ is the $(p,p')$th
  element of the matrix $B=A^2=\Sigma^{-1}$.
Similarly 
\begin{equation}
  \label{eq:tight_2}
   \sum_{k'=1}^K \{x_{n,k'}^2(t_2)-x_{n,k'}^2(t)\} =
\sum_{q,q'=1}^K b_{q,q'}
\left[\hat{V}_{n,q}(\floor{n t_2})-\hat{V}_{n,q}(\floor{n t})\right]
\left[\hat{V}_{n,q'}(\floor{n t_2})+\hat{V}_{n,q'}(\floor{n t})\right].
\end{equation}
Using the notations $\ell=\floor{n t}$, $\ell_1=\floor{n t_1}$ and $\ell_2=\floor{n t_2}$, and
decomposing the interval $[1,n]$ into $[1,\ell_1]$,
$[\ell_1+1,\ell]$, $[\ell+1,\ell_2]$ and $[\ell_2+1,n]$, we get
\begin{multline}
  \label{eq:DecompoVhat}
  \hat{V}_{n,p}(\ell)-\hat{V}_{n,p}(\ell_1) = 
  \frac{\ell-\ell_1}{n^{3/2}}\left(\sum_{i=1}^{\ell_1}h_{1,p}(X_{i,p}) \right)
  -\frac{n-(\ell-\ell_1)}{n^{3/2}}\left(\sum_{i=\ell_1+1}^{\ell}h_{1,p}(X_{i,p})\right)\\
  +\frac{\ell-\ell_1}{n^{3/2}}\left(\sum_{i=\ell+1}^{\ell_2}h_{1,p}(X_{i,p})\right)
  +\frac{\ell-\ell_1}{n^{3/2}}\left(\sum_{i=\ell_2+1}^{n}h_{1,p}(X_{i,p})\right)\;,    
\end{multline}
and
\begin{multline}
  \label{eq:DecompoVhat2}
  \hat{V}_{n,p}(\ell)+\hat{V}_{n,p}(\ell_1) = 
  \frac{\ell+\ell_1}{n^{3/2}}\left(\sum_{i=1}^{\ell_1}h_{1,p}(X_{i,p}) \right)
  -\frac{n-(\ell+\ell_1)}{n^{3/2}}\left(\sum_{i=\ell_1+1}^{\ell}h_{1,p}(X_{i,p})\right)\\
  +\frac{\ell+\ell_1}{n^{3/2}}\left(\sum_{i=\ell+1}^{\ell_2}h_{1,p}(X_{i,p})\right)
  +\frac{\ell+\ell_1}{n^{3/2}}\left(\sum_{i=\ell_2+1}^{n}h_{1,p}(X_{i,p})\right)\;,
\end{multline}
with similar results for the terms of \eqref{eq:tight_2}.
% Using the independence of the $(X_{i,k})_{1\leq i\leq n}$ and \eqref{eq:variance_h},
% \begin{equation}
%   \label{eq:majoEsperanceCarreH}
%    \PE
%   \left[
%     \left|\sum_{i=\ell_1+1}^{\ell}h_{1,p}(X_{i,p})\right|^2
%     \right]
%     \leq
%     \sum_{i=\ell_1+1}^{\ell}\PE[h_{1,p}(X_{i,p})^2]
%     =
%     (\ell-\ell_1)/3\;.
% \end{equation}
%
Equation \eqref{eq:tightnessCondition_bis} is the expected value of the product
of the squares of \eqref{eq:tight_1} and \eqref{eq:tight_2}.
Using Cauchy-Schwarz inequality, (\ref{eq:tightnessCondition_bis}) is bounded above by the sum of
several terms, obtained by inserting \eqref{eq:DecompoVhat} and
\eqref{eq:DecompoVhat2} in \eqref{eq:tight_1} and \eqref{eq:tight_2},
respectively. Among these terms, we consider the case of:
\begin{multline}
  \label{eq:majoration1}
  C_1 \sum_{p,p'=1}^K \sum_{q,q'=1}^K b_{p,p'}^2b_{q,q'}^2 
  \frac{(n-(\ell-\ell_1))^2 (\ell+\ell_1)^2 (n-(\ell_2-\ell))^2 (\ell_2+\ell)^2}{n^{12}}\\\
  \times \PE
  \left[
    \left|\sum_{i=\ell_1+1}^{\ell}h_{1,p}(X_{i,p})\right|^2
    \left|\sum_{i=1}^{\ell_1}h_{1,p}(X_{i,p}) \right|^2
    \left|\sum_{i=\ell+1}^{\ell_2}h_{1,p}(X_{i,p})\right|^2
    \left|\sum_{i=\ell_2+1}^{n}h_{1,p}(X_{i,p})\right|^2
  \right]\;.
\end{multline}
Using the independence of $(X_{i,k})_{1\leq i\leq n}$, the expected value in \eqref{eq:majoration1} can
be separated into the product of four expected values, and thus
can be bounded by 
\begin{equation}
(\ell-\ell_1)\ell_1(\ell_2-\ell)(n-\ell_2)/3^4 \leq n^2(\ell_2-\ell_1)^2/3^4.\label{eq:boundarySingleTerm}
\end{equation}
Equation \eqref{eq:majoration1} is thus bounded above by a quantity proportional
to $(\ell_2-\ell_1)^2/n^2 = (\floor{nt_2}-\floor{nt_1})^2/n^2$. All
the terms appearing in the expansion of
\eqref{eq:tightnessCondition_bis} can be treated similarly. This completes the proof of
(\ref{eq:tightnessCondition}) and thus ensures that
\begin{equation}\label{eq:conv:proc:change}
\sup_{0<t<1}\mathbf{V}_{n}(\lfloor n t\rfloor)'\Sigma^{-1}\mathbf{V}_{n}(\lfloor n t\rfloor)
  \inlaw \sup_{0<t<1}\sum_{k=1}^{K} B_{k}^{2}(t),\; n\to\infty\;.
\end{equation}
In order to prove (\ref{eq:conv:proc:change}) when $\Sigma$ is
replaced by $\hat{\Sigma}_n$, it enough to prove that 
$\sup_{0<t<1}|\mathbf{V}_{n}(\lfloor n
t\rfloor)'(\Sigma^{-1}-\hat\Sigma_n^{-1})\mathbf{V}_{n}(\lfloor n
t\rfloor)|=o_p(1).$
Note that  
$$
|\mathbf{V}_{n}(\lfloor n
t\rfloor)'(\Sigma^{-1}-\hat\Sigma_n^{-1})\mathbf{V}_{n}(\lfloor n
t\rfloor)|\leq \|\hat\Sigma_n^{-1}-\Sigma^{-1}\| \sup_{0<t<1}\|\mathbf{V}_{n}(\lfloor n
t\rfloor)\|^2\;,
$$
where $\hat\Sigma_n^{-1}\inproba\Sigma^{-1}$ and $\sup_{0<t<1}\|\mathbf{V}_{n}(\lfloor n
t\rfloor)\|^2=O_p(1)$, by (\ref{eq:conv:proc:change}), which concludes
the proof.

\bibliographystyle{chicago}
\bibliography{multirank}

\end{document}